\documentclass{article}
\usepackage{amsmath}
\usepackage{amsthm}
\usepackage{latexsym}
\usepackage{geometry}
\usepackage{amssymb}
\usepackage{cases}
\usepackage{tikz}
\title{Some Ideas on Categories and Sheaves}
\author{Dezhao Zhang}
\begin{document}
\maketitle
\begin{abstract}
We firstly introduce some key concepts in category theory, such as quotient category, completion of limits, $\mathrm{Mor}$ category, and so on; then give the concept of topology algebras and sheaves, and discuss how to restore the structue of sheaves from their stalks; lastly, we introduce the sheaf-theoretical expression for topological spaces, and rediscribe some essential items in topology and geometry by defining a kind of generally existing category sheaves.
\end{abstract}
\tableofcontents

\section{Introduction}
\paragraph{}If readers are familiar with the concepts of sheaves, (differential) manifolds, schemes and bundles, then we will naturally find the similarity between these structures, that is, they all introduce some kind of concept about gluing, and this article is making a beneficial attempt to explain this phenomenon. To achieve the goal, we firstly develop the concepts of quotient categories, completion of (co)limits, the sheaf-theoretical definition of topological spaces and so on.

\section{Basic Knowledge:Categories, Functors, Natural Transformations and Limits}
\paragraph{}To discuss the mostly used categories in ordinary math, we need push sets forward to the concept of classes.
\paragraph{}To avoid the appearance of Russell's paradox, ZFC(Zermelo-Fraenkel-Skolem)axiomatic set theory system adopts the axiom of comprehension, which claims that the method of defining sets by any first-order logic languege must rely on existing sets. In this condition, we cannot define mathematical object as "set of all sets" . Luckily, we can use axioms of classes to discribe these objects, such as the true class(class which is not set) of all sets. The behavioral pattern of classes is much similar to sets, so readers can naturally follow ways in set theory to handle classes, and use set-theretical terms and notations.
\paragraph{Definition 2.1}Let $\mathrm{Ob}$ a class, and $\{\mathrm{Mor}_{A,B}\}_{A,B\in \mathrm{Ob}}$ a family of classes  which use $\mathrm{Ob}\times \mathrm{Ob}$ as index，then we call $\mathcal{C}=(\mathrm{Ob},\{\mathrm{Mor}_{A,B}\}_{A,B\in \mathrm{Ob}})$ a graph. If for all $A,B,C\in \mathrm{Ob}$, there exists a binary function $comp:\mathrm{Mor}_{A,B}\times \mathrm{Mor}_{B,C}\to \mathrm{Mor}_{A,C}$ satisfies:
\subparagraph{1.}for all $A\in \mathrm{Ob}$, there exists identity $id_A\in \mathrm{Mor}_{A,A}$, namely $\forall f\in \mathrm{Mor}_{A,B},\mathrm{comp}(id_A,f)=f;\forall g\in \mathrm{Mor}_{B,A},\mathrm{comp}(g,id_A)=g$;
\subparagraph{2.} associative law, namely $\mathrm{comp}(\mathrm{comp}(f,g),h)=\mathrm{comp}(f,\mathrm{comp}(g,h))$;
\\\\then we call this graph a \textbf{category}.
\\\\If $\mathcal{C}=(\mathrm{Ob},\{\mathrm{Mor}_{A,B}\}_{A,B\in \mathrm{Ob}})$ is a category, then $\mathrm{Ob}$ and $\mathrm{Mor}=\bigcup_{A,B\in \mathrm{Ob}}\mathrm{Mor}_{A,B}$ can be recorded as $\mathrm{Ob}\mathcal{C}$ and $\mathrm{Mor}\mathcal{C}$, which are respectively called the object class and morphism class of $\mathcal{C}$. $\mathrm{Mor}_{A,B}$ can be recorded as $\mathrm{Mor}_{\mathcal{C}}[A,B]$, $\mathcal{C}[A,B]$ or $Hom_\mathcal{C}[A,B]$, which is called the morphism class from $A$ to $B$. Instinctively speaking, categories can be regarded as a kind of structure composed of some "points" and some "arrows" between them.
\\\\Notice that, there exist two functions $\mathrm{dom},\mathrm{cod}:\mathrm{Mor}\mathcal{C}\to \mathrm{Ob}\mathcal{C}$, if $f\in Hom_\mathcal{C}[A,B]$, then $\mathrm{dom}(f)=A,\mathrm{cod}(f)=B$, $A$ and $B$ are respectively the domain and codomain of $f$. $\mathrm{comp}$ is called the composition of morphisms, to be convenient, we mark $\mathrm{comp}(f,g)$ as $f\circ g$. So a category can be equivalently written as six-truples $(\mathrm{Ob},\mathrm{Mor},\mathrm{dom},\mathrm{cod},\mathrm{comp},\mathrm{id})$.
\\\\For convenience, in the condition of not causing ambiguity or not emphasizing morphisms, $A\in \mathrm{\mathrm{Ob}}\mathcal{C}$ is recorded as $A\in\mathcal{C}$.
\\\\A number of set-theoretical structure has its category, such as $Set$, $Grp$, $Rng$, $Top$, $K-Vect$, $Mon$, $Diff$, and so on.
\paragraph{Notice}We tacitly acknowledge the existence of identity in all pictures below.
\paragraph{Definition 2.2(special categories)}
\subparagraph{1.the empty category}whose object class and morphism class is empty set. 
\subparagraph{2.small categories}whose object class and morphism class is a set.
\subparagraph{2'.locally small categories}each of whose hom-class is a set.
\subparagraph{3.discrete categories}whose morphisms are only identity.
\subparagraph{4.single category}a category which has only one object $A$ and one morphism $id_A$, which can be recorded as $(A,id_A)$.
\subparagraph{5.simple categories}each of whose hom-class has one morphism at most.
\paragraph{Definition 2.3}Let $\mathcal{C}$ a category. We call subgraph $\mathcal{D}$ of $\mathcal{C}$ a commutiative diagram, if all paths in $\mathcal{D}$(that is a morphism-chain $(f_0,...,f_n)$ which satisfies $\mathrm{cod}(f_i)=\mathrm{dom}(f_{i+1})$) with the same start and end objects are composed to the same morphism.
\paragraph{}As for operations on categories such as subcategories $\mathcal{D}\subseteq\mathcal{C}$, product categories $\mathcal{C}\times\mathcal{D}$, coproduct categories $\mathcal{C}\sqcup\mathcal{D}$ or $\mathcal{C}+\mathcal{D}$, and opposite or dual categories $\mathcal{C}^{op}$, we don't specially introduce them in this article.
\\\\When talking about some issues in category theory, we usually set up a category as background or domain of discourse. If $P(\mathcal{C})$ is a concept/proposition/property about the domain of discourse $\mathcal{C}$, then the concept/proposition/property $P(\mathcal{C}^{op})$ is called the dual concept/proposition/property of $P(\mathcal{C})$, written as $P^{op}(\mathcal{C})$.
\paragraph{Definition 2.4(special object and morphism)} 
\subparagraph{1.}A monomorphism(or mono) is left-cancellative morphism, an epimorphism(or epic) is a right-cancellative morphism, which are dual concepts; a bimorphism is both a mono and an epic.
\subparagraph{2.}A isomorphism(or iso) has an inverse.
\subparagraph{3.}The object is called initial object that if each hom-class from all objects to it has only one morphism, the object is called terminal object that if each hom-class from all objects to it has only one morphism, which are dual concepts; the null object is both an initial and terminal object.
\paragraph{Proposition 2.1}
\subparagraph{1.}An iso has only one inverse.
\subparagraph{2.}The initial and terminal object is unique in the sense of isomorohism.
\paragraph{Definition 2.5}Let $\mathcal{C},\mathcal{D}$ are categories. If pair of functions $F=(F^{\mathrm{Ob}}:\mathrm{Ob}\mathcal{C}\to \mathrm{Ob}\mathcal{D},F^{\mathrm{Mor}}:\mathrm{Mor}\mathcal{C}\to \mathrm{Mor}\mathcal{D})$ satisfies
\subparagraph{1.}$\mathrm{dom}(F^{\mathrm{Mor}}(f))=F^{\mathrm{Ob}}(\mathrm{dom}f),\mathrm{cod}(F^{\mathrm{Mor}}(f))=F^{\mathrm{Ob}}(\mathrm{cod}f)$;
\subparagraph{2.}$F^{\mathrm{Mor}}(id_A)=id_{F^{\mathrm{Ob}}(A)}$;
\subparagraph{3.}$F^{\mathrm{Mor}}(f\circ g)=F^{\mathrm{Mor}}(f)\circ F^{\mathrm{Mor}}(g)$;
\\\\that is, preserve dom, cod, id and comp, then we call it a (convariant) functor from $\mathcal{C}$ to $\mathcal{D}$, recorded as $F:\mathcal{C}\to\mathcal{D}$. And a contravariant functor from $\mathcal{C}$ to $\mathcal{D}$ is a convariant functor from $\mathcal{C}^{op}$ to $\mathcal{D}$.
\\\\For convenience, we directly mark $F^{\mathrm{Ob}}(A)$ and $F^{\mathrm{Mor}}(f)$ as $F(A)$ and $F(f)$. Note that the function $F^\mathrm{Mor}:\mathrm{Mor}\mathcal{C}\to\mathrm{Mor}\mathcal{D}$ can be decomposed to functions $F^\mathrm{Mor}|_{[A,B]}:\mathcal{C}[A,B]\to\mathcal{D}[F(A),F(B)]$.
\paragraph{Definition 2.5'}Let $\mathcal{C}$ a category. If $F^{\mathrm{Ob}}$ and $F^{\mathrm{Mor}}$ is identity function of $\mathrm{Ob}\mathcal{C}$ and $\mathrm{Mor}\mathcal{C}$, then obviously, $F=(F^{\mathrm{Ob}}$,$F^{\mathrm{Mor}})$ is a functor, called the identity functor $id_\mathcal{C}$ of $\mathcal{C}$.
\paragraph{Definition 2.5''}Let $\mathcal{C},\mathcal{D},\mathcal{E}$ categories, $F:\mathcal{C}\to\mathcal{D},G:\mathcal{D}\to\mathcal{E}$ functors. Then obviously, pair of functions $G\circ F:=(G^\mathrm{Ob}\circ F^\mathrm{Ob},G^\mathrm{Mor}\circ F^\mathrm{Mor})$ is a functor, called the composition of $F$ and $G$. 
\\\\Now we see, functors act similarly as morphisms. In fact, use small categories as objects, and functores between them as morphisms, we get a category $Cat$, called the small catrgory category. Meanwhile, If we cancel the requirement "small", we get a "category" $CAT$, merely its object and morphism "class" are 2-classes, however its behavior is still similar, so we call it the big category category.
\paragraph{Definition 2.6(special functors)}Let $F:\mathcal{C}\to\mathcal{D}$ a functor, it is
\subparagraph{1.}faithful if each $F^{Mor}|_{[A,B]}$ is an injection;
\subparagraph{2.}injective or an embedding if is faithful and $F^{Ob}$ is an injection;
\subparagraph{3.}full if each $F^{Mor}|_{[A,B]}$ is a surjection;
\subparagraph{4.}dense if $F^{Ob}$ is a surjection in the sense of iso, that is $\forall B\in\mathcal{C},\exists A\in\mathcal{C}$ s.t. $B\cong F(A)$.
\subparagraph{5.}surjective or a projection if is full and $F^{Ob}$ is a surjection.
\paragraph{Attention}We will adopt the form of anonymous function as $F:=\begin{cases}A\mapsto F(A)\\f\mapsto F(f)\end{cases}$ to discribe functor $F$.
\paragraph{Definition 2.7}Let $F,G:\mathcal{C}\to\mathcal{D}$ functors. If function $\alpha:\mathrm{\mathrm{Ob}}\mathcal{C}\to\mathrm{Mor}\mathcal{D}$ makes $G(f)\circ\alpha(A)=\alpha(B)\circ F(f)$, then we call it a natural transformation (or nat) from $F$ to $G$, recorded as $\alpha:F\to G$.
\\\\Attention, a simple function $\mathrm{Ob}\mathcal{C}\to \mathrm{Mor}\mathcal{D}$ may be natural transformations between different functors.
\paragraph{Definition 2.7'}Let $F:\mathcal{C}\to\mathcal{D}$ a functor. If $\alpha:\mathrm{Ob}\mathcal{C}\to \mathrm{Mor}\mathcal{D}:=A\mapsto id_{F(A)}$, then obviously, it is a nat, called the identity nat $id_F$ of $F$.
\paragraph{Defination 2.7''}Let $F,G,H:\mathcal{C}\to\mathcal{D}$ functors, $\alpha:F\to G,\beta:G\to H$ nats. Then obviously, function $\beta\circ\alpha:=A\mapsto \beta(A)\circ\alpha(A)$ is a nat from $F$ to $H$, called the composition, $\circ$-product or horizontal product of $\alpha$ and $beta$.
\\\\Now we see, nats act similarly as morphisms. In fact, use functors from $\mathcal{C}$ to $\mathcal{D}$ as objects, and nats between them as morphisms, we get a catogory $Funct(\mathcal{C},\mathcal{D})$ or $\mathcal{D}^\mathcal{C}$. The existence of category $CAT$ and $Funct(\mathcal{C},\mathcal{D})$ remind us that any conclusion and concepts in category theory can be used in them, such as iso functor and iso nat (or natural iso), we don't construct them here.
\paragraph{Proposition 2.2} Union and intersection of subcategories is a subcategory. A functor's restriction on a subcategory is a functor.
\paragraph{}To discuss limits, for given background $\mathcal{C}$, we firstly introduce the diagonal functor $\Delta_\mathcal{J}:\mathcal{C}\to\mathcal{C}^\mathcal{J}$, where $\Delta_\mathcal{J}A$ is the only functor from  $\mathcal{J}$ to $(A,id_A)$, nat $\Delta_\mathcal{J}(f:A\to B):\Delta_\mathcal{J}A\to\Delta_\mathcal{J}B$ making that for all $j\in\mathcal{J}$ , it has $(\Delta_\mathcal{J}f)(j)=f$. When $\mathcal{J}$ is clear, we directly write it as $\Delta$.
\paragraph{Definition $\bf{2.8}^{\mathrm{Ob}}$}Let $\mathcal{J},\mathcal{C}\in CAT$, functor $F:\mathcal{J}\to\mathcal{C}$. Then a cone on $F$ is a two-tuples $(K\in \mathrm{Ob}\mathcal{C},\epsilon:\Delta K\to F)$, $K$ is called cone vertex.
\paragraph{Definition $\bf{2.8}^{\mathrm{Mor}}$}Let $(K,\epsilon)$ and $(K',\epsilon')$ $F$-cones, a cone morphism from $(K,\epsilon)$ to $(K',\epsilon')$ is a nat with the form as $\Delta f:\Delta K\to\Delta K'$ which makes $\epsilon'\circ\Delta f=\epsilon$.
\\\\Easy to know, cone morphism $\varphi:\Delta K'\to\Delta K$ gives a morphism $\varphi_{*}\in Hom_\mathcal{C}(K',K)$ making that $\varphi=\Delta(\varphi_{*})$; conversly, morphism $f\in Hom_\mathcal{C}(K',K)$ induces nat $\Delta f:\Delta K\to\Delta K$, and a cone $(K',\epsilon\circ \Delta f)$, and $\Delta f$ is a cone morphism between them.
\\\\Attention, the same nat can be cone morphisms between different cones. Easy to prove that all $F$-cone and cone morphisms between them form a category $con^F$, called the $F$-cone category. cocone and cocone morphism is dual concept of cone and cone morphism in the background $\mathcal{C}^\mathcal{J}$, $F$-cocone category is written as $con_F$.
\paragraph{Definition 2.9}The limit of $F$ is an $F$-con $(L,\delta)$, which makes that for each $F$-cone $(K,\epsilon)$, there exists sole con morphism $\Delta f:\Delta K\to\Delta L$, we write $(L,\delta)$ as $\lim\limits _{\longleftarrow}F$. Colimit is the dual concept of limit, recorded as $\lim\limits _{\longrightarrow} F$.
\\\\We see that the limit of $F$ is merely the terminal object in $con^F$, and its colimit is just the initial object in $con_F$. The only cone morphism from $F$-cone $(K,\epsilon)$ to $\lim\limits _{\longleftarrow}F=(L,\delta)$ (that is $\Delta f$ in Definition 2.10) is recorded as $\lim\limits _{\longleftarrow}^d\epsilon$. In the condition that not causing ambiguity or not emphasizing nats, we may call the limit its vertex as well. Let $F,G:\mathcal{J}\mathcal{C}$ have limits $(L,\delta)$ and $(L'.\delta')$, and a nat $\alpha:F\to G$, then there is naturally a morphism $\lim\limits _{\longleftarrow}\alpha:L\to L':=\lim\limits _{\longleftarrow}^d (\alpha\circ\delta)$, so we get a partial functor $\lim\limits _{\longleftarrow}:\mathcal{C}^\mathcal{J}\rightsquigarrow\mathcal{C}$, if $\mathcal{C}$ is $\mathcal{J}$-complete, then $\lim\limits _{\longleftarrow}$ is the right adjoint of $\Delta$.
\paragraph{Theorem 2.1(properties of limits)}Suppose the limit of $F$ exists and is $(L,\delta)$.
\subparagraph{1.}$\lim\limits _{\longleftarrow}F$ is unique in the sense of isomorphicness, that is, if $(L,\delta)$ and $(L',\delta')$ are both limits of $F$, then there exists iso $f:L\to  L'$ making that $\Delta f:(L,\delta)\to(L',\delta')$.
\subparagraph{2.}$\{\delta(A)\}_{A\in\mathcal{J}}$ are global monomorphisms, that is, for any two morphisms $f,g$ which $\mathrm{cod}f=\mathrm{cod}g=L$, if for all $A\in\mathcal{J}$, $\delta(A)\circ f=\delta(A)\circ g$, then $f=g$.
\subparagraph{3.}There exists a one-to-one correspondence between $F$-cones and the morphisms whose codomains are $L$. Theorem 4.2.1 is a more general result.
\\\\The proof is easy. These results have dual version of colimits.
\paragraph{Example 2.1 (special limits)}
\subparagraph{1.}Let $\mathcal{J}$ a discrete small category, the limit of $F:\mathcal{J}\to\mathcal{C}$ is $(L,\delta)$, then we call $L$ as the product of $\{F(A)\}_{A\in\mathcal{J}}$, recorded as $\prod\limits _{A\in\mathcal{J}}F(A)$, and $\{\delta(A)\}_{A\in\mathcal{J}}$ are called projection, recorded as $p_A$ or $\pi_A$; the colimit of $F:\mathcal{J}\to\mathcal{C}$ is $(L,\delta)$, then we call $L$ as the coproduct or sum of $\{F(A)\}_{A\in\mathcal{J}}$, recorded as $\coprod\limits _{A\in\mathcal{J}}F(A)$, and $\{\delta(A)\}_{A\in\mathcal{J}}$ are called embedding, recorded as $i_A$ or $\mu_A$.
\subparagraph{2.}Let $\mathcal{J}$ the category in left side below, the limit of $F:\mathcal{J}\to\mathcal{C}$ is $(L,\delta)$, then we call $\delta(A)$ as the equalizer of $F(f_1)$ and $F(f_2)$, recorded as $equ(F(f_1),F(f_2))$;  the colimit of $F:\mathcal{J}\to\mathcal{C}$ is $(L,\delta)$, then we call $\delta(B)$ as the coequalizer of $F(f_1)$ and $F(f_2)$, recorded as $coequ(F(f_1),F(f_2))$.
\[\begin{tikzpicture}
\path(0,0)node(A){$A$}(2,0)node(B){$B$}(7,0)node(a){$F(A)$}(9.5,0)node(b){$F(B)$}(5,0)node(L){$L$};
\draw[->](A)..controls(1,0.5)..node[above]{$f_1$}(B);
\draw[->](A)..controls(1,-0.5)..node[below]{$f_2$}(B);
\draw[->](a)..controls(8.25,0.5)..node[above]{$F(f_1)$}(b);
\draw[->](a)..controls(8.25,-0.5)..node[below]{$F(f_2)$}(b);
\draw[->](L)--node[above]{$\delta(A)$}(a);
\end{tikzpicture}\]
\\\\Attention that we don't particularly emphasize the morphism $\delta(A)$, because it can be derived from other mophismes.
\subparagraph{3.}Let $\mathcal{J}$ the wedge-shaped category in left side below. he limit of $F:\mathcal{J}\to\mathcal{C}$ is $(L,\delta)$, then we call $(\delta(B_1),\delta(B_2))$ as the pullback of $(F(f_1),F(f_2))$; the colimit of  $F:\mathcal{J}^{op}\to\mathcal{C}$ is $(L,\delta)$, then we call $(\delta(B_1^{op}),\delta(B_2^{op}))$ as the pushout of $(F(f_1^{op}),F(f_2^{op}))$.
\[\begin{tikzpicture}
\path (0,0)node(A){$A$}(1.5,1)node(B1){$B_1$}(1.5,-1)node(B2){$B_2$}
(4,0)node(a){$F(A)$}(6.5,1)node(b1){$F(B_1)$}(6.5,-1)node(b2){$F(B_2)$}(9,0)node(L){$L$};
\draw[->](B1)--node[above]{$f_1$}(A);
\draw[->](B2)--node[below]{$f_2$}(A);
\draw[->](b1)--node[above]{$F(f_1)$}(a);
\draw[->](b2)--node[below]{$F(f_2)$}(a);
\draw[->](L)--node[above]{$\delta(B_1)$}(b1);
\draw[->](L)--node[below]{$\delta(B_2)$}(b2);
\end{tikzpicture}\]
\\\\Attention that we don't particularly emphasize the morphism $\delta(A)$, because it can be derived from other mophismes.
\subparagraph{4.}Let us update $\mathcal{J}$ in 3. to the category in left side below (we call it multi-wedge-shaped category), there are two types of objects in $\mathcal{J}$: $\{B_\alpha\}$ and $\{A_{\alpha\beta}\}$, satisfy $A_{\alpha\beta}=A_{\alpha\beta}$ and $A_{\alpha\alpha}=B_\alpha$, and for any $\alpha,\beta$, there exists unique $f_{\alpha\beta}:B_\alpha\to A_{\alpha\beta}$, naturally $f_{\alpha\alpha}=id_{B_\alpha}$. Then the corresponding concepts are updated to paired-pushout.
\[\begin{tikzpicture}
\path (0,0)node(A){$A_{\alpha\beta}$}(2,1)node(B1){$B_\alpha$}(2,-1)node(B2){$B_\beta$}
(5,0)node(a){$F(A_{\{\alpha,\beta\}})$}(8,1)node(b1){$F(B_\alpha)$}(8,-1)node(b2){$F(B_\beta)$}(10.5,0)node(L){$L$};
\draw[->](B1)--node[above]{$f_{\alpha\beta}$}(A);
\draw[->](B2)--node[below]{$f_{\beta\alpha}$}(A);
\draw[->](b1)--node[above]{$F(f_{\alpha\beta})$}(a);
\draw[->](b2)--node[below]{$F(f_{\beta\alpha})$}(a);
\draw[->](L)--node[above]{$\delta(B_\alpha)$}(b1);
\draw[->](L)--node[below]{$\delta(B_\beta)$}(b2);
\end{tikzpicture}\]
\\\\Attention that we don't particularly emphasize the morphism $\delta(A_{\alpha\beta})$, because it can be derived from other mophismes, in addition pullback and pushout are not dual concepts.
\\\\A common situation is that $F:\mathcal{J}\to\mathcal{C}$ is an embedding, then we can regard $\mathcal{J}$ as a subcategory of $\mathcal{C}$, it's a convenient viewpoint.
\paragraph{Example 2.2(limits in concrete categories)}
\subparagraph{1.}The product in $Set$ is just Cartesian product,and the coproduct is disjoint union; the product in $Grp$ is direct product, and coproduct is free product; the product in $R-Mod$ is direct product,and the coproduct is direct sum; and so on.
\subparagraph{2.}Let $Set$ be the background. In the picture as below, $i_1,i_2,i_3,i_4$ are all inclusion maps, then $(i_1,i_2)=pullback(i_3,i_4)$ and $(i_3,i_4)=pushout(i_1,i_2)$. In fact, if we change $A\cup B$ to a set $C\supseteq A\cup B$, then $(i_1,i_2)=pullback(i_3,i_4)$ still applies.
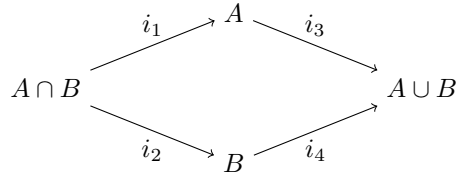
\[\begin{tikzpicture}
\path(0,0)node(a){$A\cap B$}(2.5,1)node(b){$A$}(2.5,-1)node(c){$B$}(5,0)node(d){$A\cup B$};
\draw[->](a)--node[above]{$i_1$}(b);
\draw[->](a)--node[below]{$i_2$}(c);
\draw[->](b)--node[above]{$i_3$}(d);
\draw[->](c)--node[below]{$i_4$}(d);
\end{tikzpicture}\] 
\paragraph{Theorem 2.2 (the first picture of limits)}Let $\mathcal{J},\mathcal{C}$ categories, functor $F:\mathcal{J}\to\mathcal{C}$. Then $\lim\limits _{\longleftarrow}F=(L,p\circ\Delta(e))$, where $e=equ\Bigl(\prod\limits _{f\in \mathrm{Mor}\mathcal{J}}p_{cod(f)},\prod\limits _{f\in \mathrm{Mor}\mathcal{J}}F(f)\circ p_{dom(f)}\Bigr)\in Hom(L,\prod\limits _{A\in\mathcal{J}} F(A))$. We have a similar conclusion for colimits.
\\\\The Theorem directly declares that any limits can be expressed by products and equalizers, instinctively speaking, it desposes limits to "object" part(products) and "morphism" part(equalizers).
\paragraph{Theorem 2.3 (the second picture of limits)}Let $\mathcal{J},\mathcal{C}$ categories, functor $F:\mathcal{J}\to\mathcal{C}$. Then $\lim\limits _{\longleftarrow}F=(L,\delta)$ if and only if
\subparagraph{1.}All hom-classes $Hom_{con^F}((K,\epsilon),(L,\delta))$ are not empty;
\subparagraph{2.}$\{\delta(A)\}_{A\in\mathcal{J}}$ are global monomorphisms.
\paragraph{}There are still much content in category theory, readers are advised to consult related books.

\section{Quotient Category and Sketch}
\subsection{Quotient Category}
\paragraph{Definition 3.1.1}Let $\sim^\mathrm{Ob}$ an equivalence relation on $\mathrm{Ob}\mathcal{C}$, and $\sim^\mathrm{Mor}$ an equivalence relation on $\mathrm{Mor}\mathcal{C}$, then we call the pair $\sim=(\sim^\mathrm{Ob},\sim^\mathrm{Mor})$ a precategorical equivelance relation if it satisfies
\subparagraph{1.(dom/cod-preservation)}$f\sim^\mathrm{Mor}g\Rightarrow \mathrm{dom}f\sim^{\mathrm{Ob}}\mathrm{dom}g\;\; and\;\;\mathrm{cod}f\sim^\mathrm{Ob}\mathrm{cod}g;$
\subparagraph{2.(id-preservation)}$A\sim^\mathrm{Ob}B\Rightarrow id_{A}\sim^\mathrm{Mor}id_{B};$
\subparagraph{3.(comp-preservation)} For all equivalence class $\mathbf{f},\mathbf{g}$, there exists a unique $\mathbf{h}$, making that for all $f\in\mathbf{f}$ and $g\in\mathbf{g}$ whose $\mathrm{cod}f=\mathrm{dom}g$, it has $g\circ f\in\mathbf{h}$.
\paragraph{}If $F:\mathcal{C}\to\mathcal{D}$ is a functor, then we can get a precategorical equivalence relation $\sim_F$:
$A\sim_F^\mathrm{Ob}B:=F(A)=F(B)$ and $f\sim_F^\mathrm{Mor}g:=F(f)=F(g)$. We can naturally regard precategorical equivalence relations as the popularization of set-theoretical equivalence relations.
\paragraph{}Easy to see one of consequences of dom/cod-preservation is the quotient graph $\mathcal{C}/\sim$, the graph's objects and morphisms are just equivalence classes of $\sim$, $\mathrm{dom'/cod'}\;\mathbf{f}:=\mathrm{dom/cod}f$, and naturally the projection $[-]:\mathcal{C}\to \mathcal{C}/\sim$. Another its consequence is $\sim^\mathrm{Mor}$ can be decomposed into equivalence realtions $\sim^\mathrm{Mor}_{[\mathbf{A},\mathbf{B}]}$ on the subsets of $\mathrm{Mor}\mathcal{C}$ which have the form as $[\mathbf{A},\mathbf{B}]=\bigcup\limits _{A\in\mathbf A, B\in\mathbf{B}}Hom_\mathcal{C}(A,B)$, where $\mathbf{A}$ and $\mathbf{B}$ are equivalence classes of $\sim^\mathrm{Ob}$, then graph-morphism class $Hom_{\mathcal{C}/\sim}(\mathbf{A},\mathbf{B})=[\mathbf{A},\mathbf{B}]/\sim^\mathrm{Mor}_{[\mathbf{A},\mathbf{B}]}$.
\paragraph{}Now we want to obtain a precategorical equivalence relation which is not induced from a given functor but can form a quotient category and a projection functor by itself, so we need new a condition about \textbf{feasibility}: let $\mathbf{f}$ and $\mathbf{g}$ is respectively any equivalence class of $\sim^\mathrm{Mor}_{[\mathbf{A},\mathbf{B}]}$ and $\sim^\mathrm{Mor}_{[\mathbf{B},\mathbf{C}]}$, then there exist $f\in\mathbf{f}$ and $g\in\mathbf{g}$ whose $\mathrm{cod}f=\mathrm{dom}g$, that is $f$ and $g$ can be composed. So from it we can define a partial binary function $\mathrm{comp}'(\mathbf{g},\mathbf{f}):=\mathbf{h}$, or briefly $\mathbf{g}\circ\mathbf{f}$, at this time the quotient graph $\mathcal{C}/\sim$ becomes a category, and we call $\sim$ \textbf{categorical}.
\paragraph{}A special situation is that $A\sim^{\mathrm{Ob}}B\Leftrightarrow A=B$ (we mark $\sim^{\mathrm{Ob}}$ as $=^{\mathrm{Ob}}$, which is an essantial notation), at this time the quotient category $\mathcal{C}/(=^{\mathrm{Ob}},\sim^{\mathrm{Mor}})$ is called a wide subcategory of $\mathcal{C}$.
\paragraph{Proposition 3.1.1} If an equivalence relation $\sim=(\sim^\mathrm{Ob},\sim^\mathrm{Mor})$ preserves dom/cod and comp, then it preserves id.
\begin{proof} For all $A$ and all $f$ whose $\mathrm{dom}f=A$, it has $f\circ id_A=f$, so for all $\mathbf{f}$, it has $\mathbf{f}\circ[id_A]=\mathbf{f}$, in a similar way, for all $\mathbf{g}$, $[id_A]\circ\mathbf{g}=\mathbf{g}$, so $[id_A]$ is an identity in $\mathcal{C}/\sim$, and identity is unique for each object (here is $\mathbf{A}$).
\end{proof}
\paragraph{Proposition 3.1.2}Let $F:\mathcal{C}\to \mathcal{D}$ match up with a categorical equivalence relation $\sim$, namely $\sim\Rightarrow\sim_F$, then there exists a unique $\rho:\mathcal{C}/\sim\to\mathcal{D}$ such that $\rho\circ [-]=F$.
\begin{proof}Just let $\rho([A]):=F(A)$ and $\rho([f]):=F(f)$.
\end{proof}

\subsection{Sketches and the Cochain Condition}
\paragraph{}We now discuss a case which occurs in various categories widespreadly.
\paragraph{}Consider a race of object equivalence relations on $\mathcal{C}$ which satisfy the strong isomorphism condition $A\sim^{\mathrm{Ob}}B\Rightarrow A\cong B$, it says all objects in the same equivalence class are isomorphic, we instinctly perceive that, there must be similarity between morphisms: if $A\sim^{\mathrm{Ob}}A'$，$B\sim^{\mathrm{Ob}}B'$, then there exists isomorphisms $\varphi:A\tilde{\to}A'$，$\phi:B\tilde{\to}B'$, in that way we have bijection $Hom(A,B)\to Hom(A',B'):=f\mapsto \phi\circ f\circ\varphi^{-1}$, but we still need a strengthened condition about these bijection for us to get a good morphism equivalence relation.
\paragraph{Definition 3.2.1}Let $\mathbf{A}$ an equivalence class of $\sim^{\mathrm{Ob}}$ which the strong isomorphism condition. We say a group of isomorphisms $\varphi_\mathbf{A}=\{\varphi_{AB}\in Iso(A,B) \}_{A,B\in\mathbf{A}}$ on $\mathbf{A}$ satisfies the cochain condition if 
\subparagraph{1.}$\varphi_{AA}=id_A,\forall A\in\mathbf{A}$;
\subparagraph{2.}$\varphi_{BA}\circ\varphi_{AB}=id_A,\forall A,B\in\mathbf{A}$;
\subparagraph{3.}$\varphi_{CA}\circ\varphi_{BC}\circ\varphi_{AB}=id_A,\forall A,B,C\in\mathbf{A}$.
\\\\We call it a cochain group for short.
\\\\For each equivalence class $\mathbf{A}$ of $\sim^{\mathrm{Ob}}$, we set up a cochain group $\varphi_\mathbf{A}$, and $\varphi=\{\varphi_\mathbf{A}\}_{\mathbf{A}\in\mathrm{Ob}\mathcal{C}/\sim^\mathrm{Ob}}$ called a choice of cochain groups ,then we can define a morphism equivalence relation:
\[f\sim^\mathrm{Mor}_\varphi g\Leftrightarrow\varphi_{\mathrm{cod}f,\mathrm{cod}g}\circ f\circ\varphi_{\mathrm{dom}g,\mathrm{dom}f}=g, \]
easy to prove $\sim_\varphi=(\sim^{\mathrm{Ob}},\sim^\mathrm{Mor}_\varphi)$ keeps $\mathrm{dom}$, $\mathrm{cod}$, $\mathrm{id}$ and $\mathrm{comp}$，so it's a categorical equivalence relation, therefore we get a quotient category $\mathcal{C}/\sim_\varphi$, and the quotient functor $[-]_\varphi:\mathcal{C}\to\mathcal{C}/\sim_\varphi$.
\paragraph{}
In fact, for each equivalence class $\mathbf{A}$, there must be cochain groups: choose a representative element $A$ of $\mathbf{A}$, and choose isomorphisms from $A$ to other elements as generators, then we get a cochain group by spanning the generators (here we use the idea of complete graph's minimum spanning tree). This observation gaurantees that we certainly can build a quotient category using equivalence relations which satisfy the strong isomorphism condition.
\paragraph{Lemma 3.2.1} For any choice of cochain groups $\varphi$, the quotient functor $[-]_\varphi$ is fully faithful, in other word, each $[-]_\varphi|_{[A,B]}$ is bijective.
\paragraph{}From each $\mathbf{A}\in\mathrm{Ob}\mathcal{C}/\sim^\mathrm{Ob}$, we fetch a representative element $A$, which is called a choice of representative elements $\chi$, and fetch all morphisms between these $A$s, then we obtain a full subcategory $\widetilde{\mathcal{C}}_\chi$ called a sketch of $\mathcal{C}$. An important thing is that for any choice $\chi$ and $\varphi$ as long as they are under the same $\sim^\mathrm{Ob}$, we can construct an isomorphism $\widetilde{\mathcal{C}}_\chi\cong\mathcal{C}/\sim_\varphi$ by mapping objects and morphismes in $\widetilde{\mathcal{C}}_\chi$ to their equivalence classes in $\mathcal{C}/\sim_\varphi$, and it's exactly an iso by Lemma 3.2.1. Another way is that, we can define a right inverse $S_{\varphi,\chi}$ of $[-]_\varphi$ as mapping an equivalence class of $\sim_\varphi$ to its representative element by $\chi$, easy to see that $\widetilde{\mathcal{C}}_\chi$ is exactly the image of $S_{\varphi,\chi}$, and from the factorization we obtain an iso. For this reason, we call $\mathcal{C}/\sim_\varphi$ a sketch quotient of $\mathcal{C}$.
Now we know that there is no real difference to choose different cochain groups, because the quotient categories they generate are isomophic, and \textbf{category theory doesn't care about distinction between isomorphic objects and also can't distinguish them}. 
\paragraph{}
There are two special cases of sketches:
\subparagraph{1.}$A\sim^{\mathrm{Ob}}B\Leftrightarrow A=B$, in this situation, the sketch is the original category itself.
\subparagraph{2.}$A\sim^{\mathrm{Ob}}B\Leftrightarrow A\cong B$, in this situation, sketches are skeletons.

\section{Several Categories and Functors}
\subsection{The Smooth Category and The Specture Category}
\paragraph{}All open sets on n-dimensional Euclidean space $\mathbb{R}^n$ with the standard topology, and smooth mappings between them compose a small category $Smo_n$, all open sets on all Euclidean spaces with the standard topology, and smooth mappings between them compose a small category $Smo$ called the \textbf{Smooth category}. Note that $Smo$ can be embedded into $Set$.
\paragraph{}There are several varietas of $Smo$: morphisms are continuous mappings instead of smooth mappings; objects are open sets on complex spaces and morphisms are holomorphic mappings; objects are open sets on half-spaces.
\subsection{$\mathrm{Mor}$-Categories, $\mathrm{Hom}$-Functors and $\mathrm{Dom}/\mathrm{Cod}$-Functors}
\paragraph{Definition 4.2.1}Let $\mathcal{C}$ a category. $\mathcal{C}$'s $\mathrm{Mor}$ category $\mathrm{Mor}(\mathcal{C})$ is as follows: its objects are morphisms of $\mathcal{C}$; Let $f,g$ are two objects of $\mathrm{Mor}(\mathcal{C})$, then a morphism from $f$ to $g$ is a pair of morphisms $(\varphi,\psi)$ of $\mathcal{C}$, making $\psi\circ f=g\circ\varphi$.
\paragraph{}In fact, by changing direction of $\varphi$ and $\psi$, we have anthor three similar categories which maybe not not much important here. Naturally, there are two functors 
\[\mathrm{Dom_\mathcal{C}}:\mathrm{Mor}(\mathcal{C})\to\mathcal{C}:=\begin{cases}f\mapsto \mathrm{dom}f\\(\varphi,\psi)\mapsto\varphi\end{cases} and\;\;\; \mathrm{Cod_\mathcal{C}}:\mathrm{Mor}(\mathcal{C})\to\mathcal{C}:=\begin{cases}f\mapsto \mathrm{cod}f\\(\varphi,\psi)\mapsto\psi\end{cases},\] 
and a functor \[\mathrm{Id_\mathcal{C}}:\mathcal{C}\to\mathrm{Mor}(\mathcal{C}):=\begin{cases}A\mapsto id_A\\f\mapsto(f,f)\end{cases},\]
and a partial multifunctor
\[\mathrm{Comp_\mathcal{C}}:\mathrm{Mor}(\mathcal{C})\times\mathrm{Mor}(\mathcal{C})\rightsquigarrow\mathrm{Mor}(\mathcal{C}):=\begin{cases}(f,g)\mapsto g\circ f\\ ((\varphi,\psi):f\to f',(\psi,\phi):g\to g')\mapsto(\varphi,\phi)\end{cases},\]
they reflect the struture of category. In fact, the functor $comp_\mathcal{C}$ reflects a semigroupoid structure on the object class of $\mathrm{Mor}(\mathcal{C})$, note that there is also a semigroupoid structure on the morphsim class of $\mathrm{Mor}(\mathcal{C})$, we can see that these two semigroupoids are "perpendicular" to each other from the following picture.

\[\begin{tikzpicture}
\path(0,0)node(1){}(1.5,0)node(2){}(3,0)node(3){}(0,-1.5)node(4){}(1.5,-1.5)node(5){}(3,-1.5)node(6){}(0,-3)node(7){}(1.5,-3)node(8){};
\draw[->](1)--node[above]{$\varphi$}(2);
\draw[->](2)--node[above]{$\varphi'$}(3);
\draw[->](4)--node[above]{$\psi$}(5);
\draw[->](5)--node[above]{$\psi'$}(6);
\draw[->](7)--node[above]{$\phi$}(8);

\draw[->](1)--node[left]{$f$}(4);
\draw[->](4)--node[left]{$f'$}(7);
\draw[->](2)--node[left]{$g$}(5);
\draw[->](5)--node[left]{$g'$}(8);
\draw[->](3)--node[left]{$h$}(6);
\end{tikzpicture}\] 

\paragraph{}In addition, functor $F:\mathcal{C}\to\mathcal{D}$ can be promoted to 
\[\mathrm{Mor}(F):\mathrm{Mor}(\mathcal{C})\to\mathrm{Mor}(\mathcal{D}):=\begin{cases}f\mapsto F(f)\\(\varphi,\psi)\mapsto(F(\varphi),F(\psi))\end{cases},\]
so we actually get a functor $\mathrm{Mor}:CAT\to CAT$.
\\\\And there is the diagonal plane functor \[DP_\mathcal{C}:\mathrm{Mor}^2(\mathcal{C})\to\mathrm{Mor}(\mathcal{C}):=\begin{cases}((\varphi,\psi):f\to g)\mapsto\psi\circ f\equiv g\circ\varphi\\((i,k),(m,n))\mapsto(i,n)\end{cases}.\]

\[\begin{tikzpicture}
\path(0,0)node(1){}(1.5,0)node(2){}(0,-1.5)node(3){}(1.5,-1.5)node(4){}(2,1)node(5){}(3.5,01)node(6){}(2,-0.5)node(7){}(3.5,-0.5)node(8){};
\draw[->](1)--node[above,pos=0.6]{$\varphi$}(2);
\draw[->](3)--node[below,pos=0.6]{$\psi$}(4);
\draw[->](1)--node[left,pos=0.6]{$f$}(3);
\draw[->](2)--node[right,pos=0.6]{$g$}(4);
\draw[->](1)--node[]{}(4);

\draw[->](1)--node[left,pos=0.6]{$i$}(5);
\draw[->](2)--node[right,pos=0.3]{$k$}(6);
\draw[->](3)--node[left,pos=0.75]{$m$}(7);
\draw[->](4)--node[right]{$n$}(8);

\draw[->](5)--node[above]{$\varphi'$}(6);
\draw[->](7)--node[above]{$\psi'$}(8);
\draw[->](5)--node[left]{$f'$}(7);
\draw[->](6)--node[right]{$g'$}(8);
\draw[->](5)--node[]{}(8);
\end{tikzpicture}\] 

\paragraph{}We can use these functors to rewrite the concept of natural transform:
\paragraph{Definition 4.2.2}Let $F,G:\mathcal{C}\to\mathcal{D}$ functors. Then a natural transform from $F$ to $G$ is a functor $\alpha:\mathcal{C}\to\mathrm{Mor}(\mathcal{D})$, making $\mathrm{Dom_\mathcal{D}}\circ\alpha=F$ and $\mathrm{Cod_\mathcal{D}}\circ\alpha=G$, in other words, $\alpha(f)=(F(f),G(f))$.
\paragraph{Proposition 4.2.1}Let $F,G:\mathcal{C}\to\mathcal{D}$ functors, $\alpha:F\to G$ a nat, $\Phi:\mathcal{E}\to\mathcal{C}$ a functor. Then $\alpha\circ\Phi:F\circ\Phi\to G\circ\Phi$ is a nat. 
\\A special situation is that $\mathcal{E}$ is a subcategory of $\mathcal{C}$, namely there is an inclusion functor $in:\mathcal{E}\hookrightarrow\mathcal{C}$ , so $\alpha\circ in:F\circ in\to G\circ in$ is a nat, namely the restriction $\alpha|_\mathcal{E}:F|_\mathcal{E}\to G|_\mathcal{E}$ is a nat.
\paragraph{Proposition 4.2.1'}Let $F,G:\mathcal{C}\to\mathcal{D}$ functors, $\alpha:F\to G$ a nat, $\Phi:\mathcal{D}\to\mathcal{E}$ a functor. Then $\mathrm{Mor}(\Phi)\circ\alpha:\Phi\circ F\to\Phi\circ G$ is a nat.
\paragraph{Definition 4.2.2'}Let $F,G,H:\mathcal{C}\to\mathcal{D}$ functors and $\alpha:F\to G,\beta:G\to H$ nats. Then the $\circ$-composition of $\alpha$ and $\beta$ is $\beta\circ\alpha:\mathcal{C}\to\mathrm{Mor}(\mathcal{D}):=\mathrm{comp}_\mathcal{C}\circ(\alpha\times\beta)$, which is a nat $F\to H$.
\paragraph{Definition 4.2.2''}Let $F,G:\mathcal{C}\to\mathcal{D}$ and $H,K:\mathcal{D}\to\mathcal{E}$ functors, $\alpha:F\to G$ and $\beta:H\to K$ nats. Then the $*$-composition of $\alpha$ and $\beta$ is $\beta *\alpha:\mathcal{C}\to\mathrm{Mor}(\mathcal{E}):=DP_\mathcal{E}\circ\mathrm{Mor}(\beta)\circ\alpha$, which is a nat $H\circ F\to K\circ G$.
\paragraph{}It reveals a kind of possibility: we can continue to build morphisms between nats (can be called 2-nat or 3-functor or 4-object) $\mathfrak{t}:\alpha\to\beta$, it's a functor $\mathfrak{t}:\mathcal{C}\to\mathrm{Mor}^2(\mathcal{D})$, making $\mathrm{Dom}_{\mathrm{Mor}(\mathcal{D})}\circ\mathfrak{t}=\alpha$ and $\mathrm{Cod}_{\mathrm{Mor}(\mathcal{D})}\circ\mathfrak{t}=\beta$. And similar to $\circ$-composition and $*$-composition, it has three compositions. Repeat this process again and again, we get a structure of $\mathbb{N}$ in category theory.
\paragraph{}Next we can naturally rewrite the concept of limits: so-called $con^F$ is a subcategory of $\mathrm{Mor}(\mathcal{C}^\mathcal{J})$, making $\mathrm{Dom}_{\mathcal{C}^\mathcal{J}}(con^F)=(F,id_F)$ and $\mathrm{Cod}_{\mathcal{C}^\mathcal{J}}(con^F)=\Delta(\mathcal{C})$, in other word,  $con^F=\mathrm{Dom}^{-1}_{\mathcal{C}^\mathcal{J}}(F,id_F)\cap\mathrm{Cod}^{-1}_{\mathcal{C}^\mathcal{J}}(\Delta(\mathcal{C}))$, and $\lim\limits _{\longleftarrow}F$ is $con^F$'s terminal object.
\paragraph{Theorem 4.2.1}If $F:\mathcal{J}\to\mathcal{C}$ has the limit $(L,\delta)$, then $con^F\cong \mathrm{Cod}^{-1}_{\mathcal{C}}(L,id_L)$; If $F:\mathcal{J}\to\mathcal{C}$ has the colimit $(L,\delta)$, then $con_F\cong \mathrm{Dom}^{-1}_{\mathcal{C}}(L,id_L)$.
\paragraph{Definition 4.2.3}Let $\mathcal{C}$ a locally small category. Functors \[Hom_\mathcal{C}(A,-):\mathcal{C}\to Set:=\begin{cases}B\mapsto Hom_\mathcal{C}(A,B)\\f:B\to B'\mapsto(g\in Hom_\mathcal{C}(A,B)\mapsto f\circ g)\end{cases},\] \[Hom_\mathcal{C}(-,A):\mathcal{C}^{op}\to Set:=\begin{cases}B\mapsto Hom_\mathcal{C}(B,A)\\f:B\to B'\mapsto(g\in Hom_\mathcal{C}(A,B')\mapsto g\circ f)\end{cases}.\]
\\\\When the background $\mathcal{C}$ is clear, $Hom_\mathcal{C}(A,-)$ is recorded as $Hom_A(-)$, $Hom_\mathcal{C}(-,A)$ is recorded as $Hom^A(-)$.
\paragraph{Theorem 4.2.2}Let $F,G:\mathcal{C}\to\mathcal{D}$ have limits, and $\alpha:F\to G$. Note that $\alpha:\mathcal{C}\to\mathrm{Mor}(\mathcal{D})$, then $\lim\limits _{\longleftarrow}\alpha=\lim\limits _{\longleftarrow}\alpha$, where the left is $\lim\limits _{\longleftarrow}:\mathcal{D}^\mathcal{C}\rightsquigarrow \mathcal{D}$, and the right is $\lim\limits _{\longleftarrow}:\mathrm{Mor}(\mathcal{D})^\mathcal{C}\rightsquigarrow \mathrm{Mor}(\mathcal{D})$. Colimits are similar. 
\paragraph{Theorem 4.2.3}$\lim\limits _{\longleftarrow}F=(L,\delta)\Longleftrightarrow\forall K,\lim\limits _{\longleftarrow}(Hom_K(-)\circ F)=(Hom_K(L),\mathrm{Mor}(Hom_K(-))\circ\delta);$\\\\$\lim\limits _{\longrightarrow}F=(L,\delta)\Longleftrightarrow\forall K,\lim\limits _{\longrightarrow}(Hom^K(-)\circ F)=(Hom^K(L),\mathrm{Mor}(Hom^K(-))\circ\delta).$

\subsection{$\mathcal{D}$-Obejects and $\mathcal{D}$-Categories}
\paragraph{Definition 4.3.1}Let $\mathcal{C}$ and $\mathcal{D}$ categories. We call an object $A$ on $\mathcal{C}$ is a $\mathcal{D}$-object, if $Hom^A(-):\mathcal{C}\to\mathcal{D}$; call $\mathcal{C}$ a $\mathcal{D}$-category, if all objects of $\mathcal{C}$ are $\mathcal{D}$-objects.
\paragraph{}Easy to know, if $A$ is an $Ab$-object of $\mathcal{C}$, then $(Hom_{\mathcal{C}}(A,A),+,\circ)$ is a ring, and pre-additive categories are $Ab$-categories. $Hom^A(-):\mathcal{C}\to Set$ means many concrete (locally small) categories we have ever constructed are all $Set$-categories, in the matter of category theory standard, it asserts that these categories are constructed by the set-theoretical language, and for a $Set$-object $A$, $(Hom(A,A),\circ)$ is a monoid. In addition, 2-categories are $Cat$-categories.

\subsection{The Module Category and $\mathcal{O}$-Module Sheaves}
\paragraph{}Firstly, we rewrite definition of module. Note that all homomorphisms of an abelian group consititute a ring, where its addition is to add values up $(f+g)(a):=f(a)+g(a)$ and its multiplication is composition.
\paragraph{Definition 4.4.1}Let $V$ an abelian group, $R$ a ring. A left $R$-module structure on $V$ is a ring homomorphism $f\in Hom_{Rng}(R,Hom_{Ab}(V,V))$.
\paragraph{Proposition 4.4.1}The definition above is equivalent to the classical definition.
\begin{proof}nb
\end{proof}
\paragraph{}By the way, a right $R$-module is a left $R^{op}$-module, where $R^{op}$ is the dual ring of $R$. 
\paragraph{}We could even discribe more generally this sort of mathematical objects.
\paragraph{Definition $\bf{4.4.2}^\mathrm{Ob}$}A module is a three-tuples $(R\in Rng,V\in Ab,f\in Hom(R,Hom(V,V)))$. If $V$ and $R$ is clear, we consider $f\in Hom(R,Hom(V,V)))$ as a module.
\paragraph{Definition $\bf{4.4.2}^\mathrm{Mor}$}A module morphism from $(R,V,f)$ to $(R',V',g)$ is a pair of morphisms $({\varphi\in Hom_{Rng}(R,R'),\psi\in Hom_{Ab}(V,V')})$, making that $\forall a\in R,g(\varphi(a))\circ\psi=\psi\circ f(a)$.
\\\\Easy to prove that modules and module morphisms constitute a category, called the \textbf{module category}, recorded as Mod.
\paragraph{}We see that the concept of the module category is very similar (but different) to that of $\mathrm{Mor}$-categories, so we can define similar (but different) \[\mathrm{MDom}:Mod\to Rng:=\begin{cases}(R,V,f)\mapsto R\\(\varphi,\psi)\mapsto\varphi\end{cases} and \;\;\; \mathrm{MCod}:Mod\to Ab:=\begin{cases}(R,V,f)\mapsto V\\(\varphi,\psi)\mapsto\psi\end{cases}.\] Note that $\mathrm{MDom}$ and $\mathrm{MCod}$ are not $\mathrm{Dom}_{Mod}$ and $\mathrm{Cod}_{Mod}$ which are defined on $\mathrm{Mor}(Mod)$.
\paragraph{}Using $\mathrm{MDom}$ we can restore classical:
\paragraph{Definition 4.4.3} Let $R$ a ring. Then $\mathrm{MDom}^{-1}(R,id_R)\subset Mod$ is called the left $R$-module category, recorded as $R-Mod$;
\\\\And (although we haven't mention the concept of sheaves yet):
\paragraph{Definition $\bf{4.3.4}^\mathrm{Ob}$}Let $\mathcal{O}:X\to Rng$ a ring sheaf. Then an $\mathcal{O}$-module sheaf is a module sheaf $F:X\to Mod$, making $\mathrm{MDom}\circ F=\mathcal{O}$.
\paragraph{Definition $\bf{4.3.4}^\mathrm{Mor}$} An $\mathcal{O}$-module sheaf morphism from $F$ to $G$ is a nat $\alpha:F\to G$, making $\mathrm{Mor}(\mathrm{MDom})\circ\alpha=\mathrm{Id}_{Rng}\circ\mathcal{O}$.
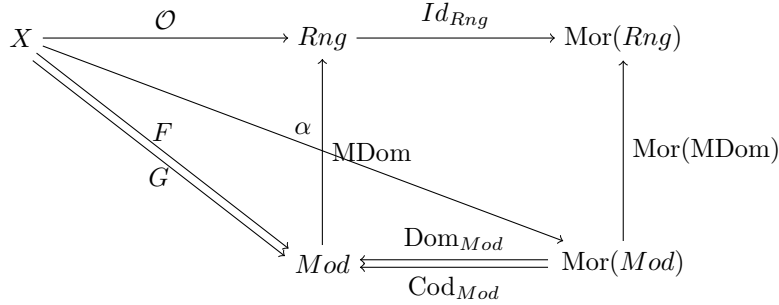
\[\begin{tikzpicture}
\path (0,0)node(X){$X$}(4,0)node(Rng){$Rng$}(8,0)node(MR){$\mathrm{Mor}(Rng)$}
(4,-3)node(Mod){$Mod$}(8,-3)node(MM){$\mathrm{Mor}(Mod)$};
\draw[->](X)--node[above]{$\mathcal{O}$}(Rng);
\draw[->](Rng)--node[above]{$Id_{Rng}$}(MR);
\draw[->](0.2,-0.2)--node[above]{$F$}(3.55,-2.8);
\draw[->](0.15,-0.3)--node[below]{$G$}(3.5,-2.9);
\draw[->](Mod)--node[right]{$\mathrm{MDom}$}(Rng);
\draw[->](X)--node[above]{$\alpha$}(MM);
\draw[->](7,-2.95)--node[above]{$\mathrm{Dom}_{Mod}$}(4.5,-2.95);
\draw[->](7,-3.05)--node[below]{$\mathrm{Cod}_{Mod}$}(4.5,-3.05);
\draw[->](MM)--node[right]{$\mathrm{Mor}(\mathrm{MDom})$}(MR);
\end{tikzpicture}\]
\paragraph{}Concepts which share the same thought and idea can be found everywhere in mathematics, such as:
\subparagraph{1.}A permutation represantation of a group, is a three-tuples \[(G\in Grp,S\in Set,f\in Hom_{Grp}(G,Aut_{Set}(S)));\]
\subparagraph{2.}A linear represantation of a group, is a three-tuples \[(G\in Grp,V\in\mathbb{R/C}-Vect,f\in Hom_{Grp}(G,Aut_{\mathbb{R/C}-Vect}(V)));\]
\subparagraph{3.}A left $R$-Algebra, is a three-tuples \[(R\in Rng,A\in Rng,f:R\to Hom_{A-Mod}(A,A)\cap Hom_{Mod-A}(A,A)).\]

\section{Topology Algebras, Sheaves and Stalks}
\subsection{Topology Algebra}
\paragraph{Definition 5.1.1}  If a small category $\mathcal{P}$ satisfies $\#Hom_{\mathcal{P}}(A,B)\bigcup Hom_{\mathcal{P}}(B,A)\leqslant 1$, then, we call it a category as partially ordered set, briefly pocategory.
\paragraph{}pocategories are simple categories, just as their name they are partially ordered sets: $A\leqslant B:=\exists f\in Hom(A,B)$, identities ensure that relation $\leqslant$ has reflexivity, condition in Definition 4.1.1 ensures antisymmetry, and composition ensures transitivity. The other way round, any poset is a pocategory: there exists a unique morphism in $Hom(A,B)$ if $A\leqslant B$. Therefore posets and pocategories are equivalent concepts, and convariant functors between pocategories are equivalent to monotonically increasing functions between posets. But by using concept of pocategories we obtain tools of category theory (especially functors and limits) to discuss them.
\paragraph{Convetion} We conveniently write the unique morphism in $Hom_P(A,B)$ as $A\leqslant B$ if exists.
\paragraph{Definition $\bf{5.1.2}^\mathrm{Ob}$}Let $X$ a set. If 0-ary operations(canstants) $1_X,0_X\in X$, 2-ary operations $\wedge:X\times X\to X$ and a subset operation $\bigvee:2^X\to X$ satisfy
\subparagraph{1.}commutative and associative law of $\wedge$;
\subparagraph{2.}absorptive law: $\bigvee\{x,x\wedge y\}=x,x\wedge\bigvee\{ x,y_\alpha\}_{\alpha\in J}=x;$
\subparagraph{3.}distributive law: $\bigvee\{x,y\wedge z\}=\bigvee\{x,y\}\wedge\bigvee\{x,z\},x\wedge\bigvee\{y_\alpha\}_{\alpha\in J}=\bigvee\{x\wedge y_\alpha\}_{\alpha\in J};$
\subparagraph{4.}identity law: $1_X\wedge x=x,\bigvee\{0_X,x\}=x;$
\\\\then we call the five-tuples $(X,1_X,0_X,\bigvee,\wedge)$ as a topology algebra.
\\\\The concept of topology algebra is connotatization of topology structure or open sets of topological space, and we will see that topological space can be restored from the topology algebra its open sets form.
\\\\Notice that from two absorption laws we can infer that $\wedge$ and $\bigvee$ are idempotent: $x\wedge x=x,\bigvee\{x,x\}\equiv\bigvee\{x\}=x$. Likewise notice that from two absorption laws and identity laws we can infer that the constants are absorbing elements: $\bigvee\{1_X,x\}=1_X,0_X\wedge x=0_X$.
\\\\To be convenient, if $\bigvee$ acts on finite set $\{x_i\}_{i=1}^n$, then we mark $\bigvee\{x_i\}_{i=1}^n$ as $x_1\vee\dots\vee x_n$, it's rational.
\\\\Likewise, on account of idempotency, commutativity and associativity of $\wedge$, it can be naturally expanded to functions on finite subsets of $X$: $\bigwedge\{x_i\}_{i=1}^n :=\bigwedge\limits _{i=1}^n x_i=x_1\wedge\dots\wedge x_n$; the other way round, $x\wedge y:=\bigwedge\{x,y\}$. They are equivlent expression.
\paragraph{Definition $\bf{5.1.2}^\mathrm{Mor}$}A topology algebra homomorphism $f:X\to Y$ is a function satisfying
\subparagraph{1.}$f(1_X)=1_Y,f(0_X)=0_Y;$
\subparagraph{2.}$f(x\wedge y)=f(x)\wedge f(y);$
\subparagraph{3.}$f(\bigvee\limits _{\alpha\in J}x_\alpha)=\bigvee\limits _{\alpha\in J}f(x_\alpha).$
\\\\Any topology algebra is a poset: $x\leqslant y:=x\vee y=y(\Leftrightarrow x\wedge y=x)$, and easy to see $1_X$ and $0_X$ is respectively maximum and minimum of $(X,\leqslant)$, $\bigvee$ and $\wedge$ is respectively supremum and finite infimum on $(X,\leqslant)$. The other way round, if a poset has supremum and finite infimum, then it's a topology algebra. They are equivalent expression. Since topology algebras are posets, they are pocategories, therefore we have a definition in category theory:
\paragraph{Definition 5.1.3}If a pocategory $X$ satisfies 
\subparagraph{1.}having the (initial) terminal object, particularly recorded as $(0_X)\;1_X$;
\subparagraph{2.}finite product complete, particularly recorded as $\wedge$;
\subparagraph{3.}coproduct complete, particularly recorded as $\bigvee$;
\subparagraph{4.}absorptive, distributive and identity law;
\\\\then we call it a category as topology algebra, briefly TA category. Naturally, morphisms are functors preserving these limits, and they compose the \textbf{category of topology algebras}, recorded as $\mathbf{TAlg}$.
\paragraph{Proposition 5.1.1 (definition and properties of subalgebra)}Let $X\in TAlg$,
\subparagraph{1.(definition)}Let $y\in X$, then subcategory $\{x\in X|x\leqslant y\}$ is still a algebra, called a subalgebra of $X$, recorded as $y^\leqslant$. The other way round, suppose $Y$ is a subcategory of $X$, then $Y=(1_Y)^\leqslant$; to sum up, $1_Y=y\Leftrightarrow\ Y=y^\leqslant$;
\subparagraph{2.(compatible property)}the intersection category $Y\cap Z$ of any two subalgebras $Y$ and $Z$ of $X$ is subalgebra, and $Y\cap Z=(1_Y\wedge 1_Z)^\leqslant$;
\subparagraph{3.(gluing property)}Let $\{Y_\alpha\}_{\alpha\in J}$ a family of subalgebras of $X$, then there exists a unique subalgebra $Z$, making that $\bigcup\limits _{\alpha\in J}Y_\alpha$ is a integrally confinal subset of  $Z=(\bigvee\limits _{\alpha\in J}1_{Y_\alpha})^\leqslant$, that is $\forall x\in Z-\{0_X\},\exists y\in\bigcup\limits _{\alpha\in J}Y_\alpha-\{0_X\},y\leqslant x$. We write $Z$ as $\bigcup\limits _{\alpha\in J}^\circ Y_\alpha$.
\paragraph{Attention}Not all subcategories which are exactly topology algebras are subalgebras.
\subsection{Sheaves and Cosheaves}
\paragraph{Definition  $\bf{5.2.1}^\mathrm{Ob}$}Let $X$ a TA category. A $\mathcal{C}$-presheaf on $X$ is a contravariant functor from $\mathcal{C}$ to $X$.
\paragraph{Definite $\bf{5.2.1}^\mathrm{Mor}$}A presheaf morphism betwwen $F,G:X\to\mathcal{C}$ is simply a nat $\alpha:F\to G$.
\\\\$\mathcal{C}$-presheaves on $X$ and morphisms compose the category $PSh^{\mathcal{C}}(X)\equiv Funct(X^{op},{\mathcal{C}})$, Similarly, $CoPSh^{\mathcal{C}}(X)\equiv Funct(X,\mathcal{C})$.
\paragraph{Definition 5.2.2}We call presheaf $F:X\to\mathcal{C}$ a sheaf, if it sytisfies
\subparagraph{1.}$F(0_X)$ is the terminal object in $\mathcal{C}$.
\subparagraph{2.(the gluing axiom)}For any $x=\bigvee x_\alpha$, it has $(x_\alpha\leqslant x)=paired-pullback(F(x_\alpha\wedge x_\beta\leqslant x_\alpha),F(x_\alpha\wedge x_\beta\leqslant x_\beta))$.
\\\\The definition has acquiesced in these completeness of $\mathcal{C}$.
\paragraph{}The definition of cosheaves is similar. Defintion 5.2.3 is called the categorical defintion of (co)sheaves, we advise readers to consult geometry books to learn about the set-theoretical definition of (co)sheaves in a concrete category, which is usually $Set$, $Ab$, and $Rng$, especially terms such as section, restriction mapping, compatible and gluing axiom.
\paragraph{}morphisms between (co)sheaves are simply nats just like (co)presheaves. The category of $\mathcal{C}$-sheaves on $X$ denotes as $Sh^\mathcal{C}(X)$, and the category of $\mathcal{C}$-cosheaves denotes as $CoSh^\mathcal{C}(X)$. We see $Sh^\mathcal{C}(X)\subset PSh^\mathcal{C}(X)\subset Funct(X^{op},\mathcal{C})$ and $CoSh^\mathcal{C}(X)\subset CoPSh^\mathcal{C}(X)\subset Funct(X,\mathcal{C})$ are chains of full subcategory.
\paragraph{Definition 5.2.3}For any algebra homomorphism $f:Y\to X$, $\mathcal{C}$-presheaves $F,G$ on $X$, and nat $\alpha:F\to G$, easy to see that $F\circ f,G\circ f$ are $\mathcal{C}$-presheaves on $Y$ and $\alpha\circ f:F\circ f\to G\circ f$ is a nat, therefore it actually provides a contravariant functor $PSh^{\mathcal{C}}(-):TAlg\to CAT$.
\\\\In a similar way, we have contravariant functors $CoPSh^{\mathcal{C}}(-)$, $Sh^{\mathcal{C}}(-)$ and $CoSh^{\mathcal{C}}(-)$. We see that the functor $PSh^{\mathcal{C}}(-)$ is similar to $Hom^A(-)$. 
\paragraph{Theorem 5.2.1 (definition and properties of local (co)presheaves)}
\subparagraph{1. (definition)}Let $F\in PSh^{\mathcal{C}}(X)$, $Y$ is a subalgebra of $X$, then $F$'s restriction $F|_Y\in PSh^{\mathcal{C}}(Y)$
on $Y$ is called a local prsheaf of $F$, it can also be written as $F|_{1_Y}$; 
\subparagraph{2. (gluing of local sheaves)}Let $\mathcal{C}$ a paired-pullback-complete category, $\{Y_\alpha\}_{\alpha\in J}$ a family of subalgebras of $X$, and $\{F_\alpha\in Sh^\mathcal{C}(Y_\alpha)\}_{\alpha\in J}$ a family of compatible sheaves, that is $\forall\alpha,\beta\in J,F_\alpha|_{Y_\alpha\cap Y_\beta}=F_\beta|_{Y_\alpha\cap Y_\beta}$. Then there exists a unique sheaf $F\in Sh^\mathcal{C}(Z)$, making $\forall\alpha\in J,F|_{Y_\alpha}=F_\alpha$, where $Z=\bigcup\limits ^\circ_{\alpha\in J} Y_\alpha$.
\paragraph{}The proposition is still true by changing $PSh$ to $CoPSh$, $CoSh$ and $Sh$ in 1; and $Sh$ to $CoSh$ in 2.
\paragraph{Definition $\bf{5.2.4}^\mathrm{Ob}$}Let $X$ a topology algebra, $F$ a $\mathcal{C}$-(co)sheaf on $X$. Then we call two-tuples $(X,F)$ a $\mathcal{C}$-(co)sheaved algebra.
\paragraph{Definition $\bf{5.2.4}^\mathrm{Mor}$}Let $(X,F)$ and $(Y,G)$ $\mathcal{C}$-(co)sheaved algebras. Then a $\mathcal{C}$-(co)sheaved algebra homomorphism from $(X,F)$ to $(Y,G)$ is a two tuples $(f,\alpha)$, where $f:Y\to X$ is a topology algebra homomorphism, and $\alpha:F\circ f\to G$ is a nat (namely (co)sheaf morphism).
\paragraph{}They compose the category of $\mathcal{C}$-(co)sheaved algebras $(Co)Sh^\mathcal{C}$. Turning the direction of $\alpha$ to $G\to F\circ f$, we get the category ${}^\mathcal{C}(Co)Sh$. Notice that $(Co)Sh^\mathcal{C}$ is similar to $Mor$-categories and $Mod$, so similarly, there are functors 
\[\mathrm{(Co)SDom^\mathcal{C}}:(Co)Sh^\mathcal{C}\to TAlg:=\begin{cases}(X,F)\mapsto X\\(f,\alpha)\mapsto f\end{cases},\] 
and 
\[\mathrm{(Co)SCod^\mathcal{C}}:(Co)Sh^\mathcal{C}\to\mathcal{C}:=\begin{cases}(X,F)\mapsto F(1_X)\\(f,\alpha)\mapsto\alpha(1_Y)\end{cases}.\]
\paragraph{}For any $X\in TAlg$, there is a sheaf\[Sub(X):X\to TAlg:=\begin{cases}x\mapsto x^\leqslant\\x_1\leqslant x_2\mapsto (x\in x_2^\leqslant\mapsto x\wedge x_1)\end{cases},\]
and for any algebra homomorphism $f:Y\to X$, there is a sheaved algebra homomorphism \[Sub(f):Sub(Y)\to Sub(X)\circ f:=y\mapsto f|_{y\leqslant},\]
so it provides a functor $Sub:TAlg\to {}^{TAlg}Sh$, which reflects natural properties of topology algebras.. 
\subsection{Particles}
\paragraph{Definition 5.3.1}Let $X$ a topology algebra. If $p$ is a subset of $X-\{0_X\}$, and satisfies
\subparagraph{1. ($\wedge$-closed or directed)}$x,y\in p\Rightarrow x\wedge y\in p$,
\subparagraph{2. (strong locally cofinal)}$\bigvee x_\alpha\in p\Rightarrow \exists x_\alpha\in p$,
\subparagraph{3. (upward-closed)} $x\in p\;\;and\;\;y>x\Rightarrow y\in p$,
\\\\then we call it a particle of $X$. All particles of $X$ denote as $Patl_X$.
\paragraph{Theorem 5.3.1}Let $x\in X\in TAlg$, the set representation of $x$ is defiend as $T^X_x:=\{p\in Patl_X|x\in p\}$ (when $X$ is clear, we can briefly write it as $T_x$), or $x\in p\Leftrightarrow p\in T_x$. Then $T_{x\wedge y}=T_x\cap T_y$, $T_{\bigvee x_\alpha}=\bigcup T_{x_\alpha}$.
\begin{proof}
\textbf{1.}From the definition we know that $p\in T_x\cap T_y$ iff $x\in p$ and $y\in p$, because $p$ is directed,  it has $x\wedge y\in p$, namely $p\in T_{x\wedge y}$, therefore $T_x\cap T_y\subseteq T_{x\wedge y}$. From Lemma 5.3.1 we know that $x\wedge y\in p\Rightarrow x\in p$, so $T_{x\wedge y}\subseteq T_x$, for the same reason it has $T_{x\wedge y}\subseteq T_y$, therefore $T_{x\wedge y}\subseteq T_x\cap T_y$. To sum up, $T_{x\wedge y}= T_x\cap T_y$.

\textbf{2.}$p\in\bigcup T_{x_\alpha}$ iff there exists $x_\alpha\in p$, from Lemma 5.3.1 we know $\bigvee x_\alpha\in p$, namely $p\in T_{\bigvee x_\alpha}$, therefore $\bigcup T_{x_\alpha}\subseteq T_{\bigvee x_\alpha}$. We can straightway get $T_{\bigvee x_\alpha}\subseteq\bigcup T_{x_\alpha}$ from the definition. To sum up, $T_{\bigvee x_\alpha}=\bigcup T_{x_\alpha}$.
\end{proof}
\paragraph{}The most important result of Theorem 5.3.1 is that we get a cosheaf
\[T^X:X\to Set:=\begin{cases}x\mapsto T^X_x\\x\leqslant y\mapsto T^X_x\subseteq T^X_y
\end{cases},\]
and a topological space $(Patl_X,\{T_x\}_{x\in X})$, which is an inherent structure of $X$.
Let $f:Y\to X$ be an algebra homomorphism, we define an continuous mapping as
\[Patl_f:Patl_X\to Patl_Y:=p\mapsto f^{-1}(p),\]
where $f^{-1}(p)$ is indeed a particle: let $x,y,x_\alpha\in Y$, then $f(x)\in p$ and $f(y)\in p\Rightarrow f(x\wedge y)=f(x)\wedge f(y)\in p$; $f(\bigvee x_\alpha)\in p\Rightarrow \bigvee f(x_\alpha)\in p\Rightarrow \exists f(x_\alpha)\in p$; $f(x)\in p$ and $y>x\Rightarrow f(y)>f(x)\in p$; and the inverse image of $T_y$ is obviously $T_{f(y)}$. So we get a contravariant functor $Patl_{-}:TAlg\to Top$, which is the adjoint of the forgotten functor $OP_{-}:Top\to TAlg$, where $OP$ means open sets. Similarly, we have
\[T^f:T^X\circ f\to T^Y:=y\mapsto(p\in T^X_{f(y)}\mapsto f^{-1}(p)\in T^Y_y),\]
it gives a functor $T^{-}:TAlg\to CoSh^{Set}$.
\paragraph{Proposition 5.3.1}Let $X$ an algebra, then $T^X_{-}:X\to (OP_{-}\circ Patl_{-})(X)$ is an isomorphism iff $T_x=T_y\Rightarrow x=y$, namely it's injective.
\paragraph{}We call such algebras topological. Now we consider another direction, let $(M,OP_M)$ a topological space, for any point $x\in M$, we define a particle on $OP_M$: $p^M_a:=\{U\in OP_M|a\in U\}$, which gives a continuous function $p^M_{-}:M\to (Patl_{-}\circ OP_{-})(M)$, we can briefly write it as $p_{-}$ when $M$ is clear.
\paragraph{Proposition 5.3.2} $p_{-}$ is a homeomorphism iff every particle of $OP_M$ has the form $p_a$ and $p_a=p_b\Rightarrow a=b$, namely it's bijective.
\\\\We call such spaces sober.
\paragraph{Proposition 5.3.3}Hausdorff implies sober.
\begin{proof}Suppose $p_{-}$ is not injective, we can choose two distinct points $a\ne b\in M$ making that $p_a=p_b=p$, so there is no $U,V\in p$, namely $a\in U$ and $b\in V$, such that $U\cap V=\emptyset$; suppose $p_{-}$ is not surjective, let $p$ be a particle which isn't in the form of $p_a$, suppose $\bigcap_{U\in p}U\ne\emptyset$, then choose an element $a$ in it, so $p\subset p_a$,   
\end{proof}
\paragraph{Proposition 5.3.4}$p_{-}$ is surjective iff it's a quotient mapping.

\subsection{Stalks}
\paragraph{Definition 5.4.1}Let $F:X\to\mathcal{C}$ a (co)presheaf, then the (co)stalk of $F$ on $p\in Patl_X$ is defined as $F_p:=(\lim\limits _{\longleftarrow})\lim\limits _{\longrightarrow} F|_p$.
\paragraph{Definition 5.4.2}Presheaf $F:X\to\mathcal{C}$ is called preapex, if for any $x\in X$, $F(x\leqslant 1_X)$ is an epic in $\mathcal{C}$. A preapex presheaf is called apex, if for any $x\ne y$, there doesn't exist an iso $f:F(x)\to F(y)$ such that $f\circ F(x\leqslant 1_X)=F(y\leqslant 1_X)$. Similarly, replacing presheaf by copresheaf and epic by mono, we get concept of apex copresheaf.
\\\\If $F:X\to\mathcal{C}$ is an apex presheaf where $\mathcal{C}$ is a concrete category, then all of its maximal sections are global sections, which means it's explicit:  all information about it is stored and only stored in its global sections $F(1_X)$, that's why we call it "apex".
\paragraph{Proposition 5.4.1}Let $F$ a preapex presheaf, then any $F(x\leqslant y)$ is an epic.
\begin{proof}Firstly we prove a simple conclusion: if $f=g\circ h$ is an epic, then $g$ is an epic. Suppose $g$ is not an epic, then there exist two different morphisms $k_1,k_2$ making that $k_1\circ g=k_2\circ g$, which means $k_1\circ f=k_1\circ g\circ h=k_2\circ g\circ h=k_2\circ f$, then $f$ is not an epic contradicting the condition. Then we can easily see it from $F(x\leqslant y)\circ F(y\leqslant 1_X)=F(x\leqslant 1_X)$.
\end{proof}
\paragraph{}To discuss relations between presheaves and sheaves, we need concepts as follows. 
\paragraph{Definition 5.4.3} Let $\mathcal{C}$ a category and $B\in\mathcal{C}$ an object. A subobject of $B$, or $B$-subobject is a pair $(A\in\mathcal{C},i:A\rightarrow B)$, where $i$ is a mono; an $B$-submorphism between $B$-subobjects $(A,i)$ and $(A',i')$ is a morphism $j:A\to A'$ such that $i=i'\circ j$.
\\\\Easy to prove that submorpshims are monos, so submorphisms provide some new subobject-reletions between subobjects. Easy to see that $B$-subobjects and $B$-submorphisms compose a category $B-Sub=\mathrm{Cod}_\mathcal{C}^{-1}(B,id_B)|_{mono}$.
\paragraph{Definition 5.4.4}Let $\{(A_\alpha,\rho_\alpha)\}$ a family of $B$-subobjects. The union of them is a $B$-subobjects $(\bigcup A_\alpha,\bigcup \rho_\alpha)$, such that for each $\alpha$, there is a $B$-submorphism from $(A_\alpha,\rho_\alpha)$ to $(\bigcup A_\alpha,\bigcup \rho_\alpha)$, and for any other $B$-subobject which satisfies this condition, there exists a unique $B$-submorphism from $(\bigcup A_\alpha,\bigcup \rho_\alpha)$ to it, in other word, the union is the "smallest" $B$-subobject which is "bigger" than all $(A_\alpha,\rho_\alpha)$; Ditto with all the $B$-submorphisms reversed, we get the intersection $(\bigcap A_\alpha,\bigcap \rho_\alpha)$, in other word, the intersection is the biggest $B$-subobject which is smaller than all $(A_\alpha,\rho_\alpha)$; The differece $(A_\alpha-A_\beta,\rho_\alpha-\rho_\beta)$ is the smallest $B$-subobject which can be unioned to $(A_\alpha\cup A_\beta,\rho_\alpha\cup\rho_\beta)$ with $(A_\beta,\rho_\beta)$.
\paragraph{Definition 5.4.5}Let $f:A\to B$ a morphism. Then an image of $f$ is a $B$-subobject $(\mathrm{Im}f,\rho)$ making that there is a morphism $\pi:A\to\mathrm{Im}f$ such that $f=\rho\circ\pi$, which is called a factorization of $f$ by $\rho$, and for any other $B$-subobject by which $f$ can be factorized, there is only one $B$-submorphism from $(\mathrm{Im}f,\rho)$ to it. In other word, the image is the smallest $B$-subobject by which $f$ can be factorized.
\paragraph{Proposition 5.4.2}The union, intersection, defference and image is unique in the sense of isomorphism, if exists. 
\paragraph{Proposition 5.4.3}If $\{(A_\alpha,i_\alpha)\}_{\alpha\in J}$ is a family of $B$-subobjects and $I\subseteq J$, then there is a unique $B$-submorphism $I\hookrightarrow J:\bigcup_{\alpha\in I}A_\alpha\to\bigcup_{\alpha\in J}A_\alpha$. There are similar conclusions for intersections, differences and images.
\paragraph{Proposition 5.4.4}Let $(A_\alpha,i_\alpha)$ a family of $B$-subobjects and $(C_\alpha,j_\alpha)$ a family of $D$-subobjects, if there exist morphisms $f_\alpha:A_\alpha\to C_\alpha$ and $f:B\to D$ such that $f\circ i_\alpha=j_\alpha=f_\alpha$, then there is a unique $\bigcup f_\alpha:\bigcup A_\alpha\to C_\alpha$, making all diagrams commutative.
\\\\In the view of Mor cateogries, $(f_\alpha,(i_\alpha,j_\alpha))$ is a subobject of $f$ in the Mor catrgory of the background category, and the theorem is saying that their union $\bigcup f_\alpha$'s domain and codomain is $\bigcup A_\alpha$ and $\bigcup B_\alpha$, in other word, the two union operations is commutative with the functors $\mathrm{Dom/Cod}$. Of course, there are similar conclusions for intersections, differences and images.
\paragraph{Definition 5.4.2}Let $F:X\to\mathcal{C}$ a (co)presheaf. The fiber space of $F$ is the copresheaf \[F_{fib}:X\to\mathcal{C}:=\begin{cases}x\mapsto\coprod\limits _{p\in T_x}F_p\\(x\leqslant y)\mapsto(i_{xy}:\coprod\limits _{p\in T_x}F_p\to\coprod\limits _{p\in T_y}F_p)\end{cases},\]
where $i_{xy}$ is embedding: due to $T_x\subseteq T_y$, $\coprod\limits _{p\in T_y}F_p=\coprod\limits _{p\in T_x}F_p+\coprod\limits _{p\in T_y-T_x}F_p$;
\\\\the section space of $F$ is the presheaf
\[F_{sec}:X\to\mathcal{C}:=\begin{cases}x\mapsto\prod\limits _{p\in T_x}F_p\\(x\leqslant y)\mapsto(p_{xy}:\prod\limits _{p\in T_y}F_p\to\prod\limits _{p\in T_x}F_p)\end{cases},\]
where $p_{xy}$ is projection: due to $T_x\subseteq T_y$, $\prod\limits _{p\in T_y}F_p=\prod\limits _{p\in T_x}F_p\times\prod\limits _{p\in T_y-T_x}F_p$.
\\\\Easy to see that for any presheaf $F:X\to\mathcal{C}$, let $(F_p,\delta_p)=\lim\limits _{\longrightarrow}F|_p$, then there is naturally a nat $\alpha^F_{sec}:F\to F_{sec}:=x\mapsto\prod\limits _{p\in T_x}\delta_p(x)$, and for copresheaf $F$, there is $\alpha^F_{fib}:F_{fib}\to F:=x\mapsto\coprod\limits _{p\in T_x}\delta_p(x)$. 
\paragraph{Theorem 5.4.1}For any (co)presheaf $F$, $F_{fib}$ is a cosheaf and $F_{sec}$ is a sheaf.
\begin{proof}let's regard $T_x$ as a discrete category, then we define a functor $G:T_x\to\mathcal{C}:=p\mapsto F_p$. Now we prove that for $A\in\mathcal{C}$, there is a $G$-cone whose vertex is $A$, if and only if, for all $x=\bigvee x_\alpha$, there is a cone of multi-wedge-shaped subcategory $\mathcal{D}=(\{F_{sec}(x_\alpha)\},\{F_{sec}(x_\alpha\wedge x_\beta)\},\{F_{sec}(x_\alpha\wedge x_\beta\leqslant x_\alpha)\})$ whose vertex is $A$. 

$\Rightarrow$ If $(A,\epsilon)$ is a $G$-cone, then $\epsilon|_{T_{x_\alpha}}$ is a $G|_{T_{x_\alpha}}$-cone, it's rational because $T_{x_\alpha}\subseteq T_x$, noitce that $\lim\limits _{\longleftarrow}G|_{T_{x_\alpha}}=F_{sec}(x_\alpha)$, so there is a unique $\epsilon'(x_\alpha):A\to F_{sec}(x_\alpha)$, $(A,\epsilon')$ is exactly a $\mathcal{D}$-cone because of some properties of product (such as commutative and associative), $T_{x_\alpha\wedge x_\beta}=T_{x_\alpha}\cap T_{x_\beta}$ and $T_x=\bigcup T_{x_\alpha}$.

$\Leftarrow$ If $(A,\epsilon)$ is a $\mathcal{D}$-cone, then we define $\epsilon'(p)=\pi_\alpha(p)\circ\epsilon(x_\alpha)$ for $p\in T_{x_\alpha}$, where $(F_{sec},\pi_\alpha)=\lim\limits _{\longleftarrow}G|_{T_{x_\alpha}}$, it is rational because $T_{x_\alpha\wedge x_\beta}=T_{x_\alpha}\cap T_{x_\beta}$ and $T_{x_\alpha}\subseteq T_x$, and $(A,\epsilon')$ is exactly a $G$-cone because $T_x=\bigcup T_{x_\alpha}$.

So now we directly get that $\lim\limits _{\longleftarrow}\mathcal{D}=\lim\limits _{\longleftarrow}G=F_{sec}(x)$ which is exactly the gluing axiom.

Notice that $T_{0_X}=\emptyset$, so $G_{0_X}:T_{0_X}\to\mathcal{C}\equiv\emptyset_\mathcal{C}$, according to Proposition 4.1.1, the limit of $G_{0_X}$, namely the $F_{sec}(0_X)$, is exactly the terminal object of $\mathcal{C}$.

The proof of $F_{fib}$ is a cosheaf is similar.
\end{proof}
\paragraph{}Easy to see that, if $\mathcal{C}$ has the property that embbedings are monos and projections are epics, then $F_{fib}$ and $F_{sec}$ are preapex, further, if $X$ is topological, then they are apex.

\paragraph{Theorem 5.4.2}If $\mathcal{C}$ is a paired-pullback complete category and has the terminal object, then for any presheaf $F\in PSh^{\mathcal{C}}(X)$, there exists a unique sheaf $\overline{F}\in Sh^{\mathcal{C}}(X)$ and a nat $\theta_F:F\to\overline{F}$, making that for any sheaf $G$ and nat $\alpha:F\to G$, there exists a unique nat $\overline{\alpha}:\overline{F}\to G$, making $\overline{\alpha}=\alpha\circ\theta_F$.
\begin{proof}For arbitrary $x\in X$, let $O=\{x_\alpha\}$ a covering of it: $x=\bigvee x_\alpha$. Consider the multi-wedge-shaped subcategory $(\{F(x_\alpha)\},\{F(x_\alpha\wedge x_\beta)\},\{F(x_\alpha\wedge x_\beta\leqslant x_\alpha)\})$, we write its limit as $(F_O(x),x_\alpha\mapsto F_O(x_\alpha\leqslant x))$. Notice that $(F(x),x_\alpha\mapsto F(x_\alpha\leqslant x))$ is a cone of that subcategory, so there is a $R_O(x):F(x)\to F_O(x)$. Notice that morphisms $\{F_O(x_\alpha\leqslant x)\circ\delta_p(x_\alpha)\}_{x_\alpha\in O, p\in T_{x_\alpha}}$ compose a cone of discrete subcategory $\{F_p\}_{p\in T_x}$ due to $\bigcup T_{x\alpha}=T_x$, so there is a morphism $S_O(x):F_O(x)\to F_{sec}(x)$, easy to see that $S_O(x)\circ R_O(x)=\alpha^F_{sec}$ for arbitrary covering $O$. Let $(\overline{F}(x)=\bigcup\mathrm{Im}\;S_O(x),\rho(x))$ be the union of $S_O(x)$'s images over all coverings of $x$, and $\pi_O(x):F_O(x)\to \overline{F}(x)$ such that $\rho(x)\circ\pi_O(x)=S_O(x)$, because $\rho(x)$ is a mono, so $\theta_F(x)=\pi_O(x)\circ R_O(x)$ is the same for arbitrary coverings $O$.

For $y\leqslant x$, easy to see that $O\wedge y=\{y_\alpha=x_\alpha\wedge y\}$ is a covering of $y$, so $\{F(x_\alpha\wedge y\leqslant x_\alpha)\circ F_O(x_\alpha\leqslant x)\}$ compose a cone of multi-wedge-shaped subcategory $(\{F(x_\alpha\wedge y)\},\{F(x_\alpha\wedge x_\beta\wedge y)\},\{F(x_\alpha\wedge x_\beta\wedge y\leqslant x_\alpha\wedge y)\})$, therefore there is an $F_O(y\leqslant x):F_O(x)\to F_{O\wedge y}(y)$, easy to prove that $F_{sec}(y\leqslant x)\circ S_O(x)=S_{O\wedge y}(y)\circ F_O(y\leqslant x)$, so there is the morphism $\mathrm{Im}(S_O(x),S_{O\wedge y}(y)):\mathrm{Im}\;S_O(x)\to\mathrm{Im}\;S_{O\wedge y}(y)$, and $\overline{F}(y\leqslant x)=(\{O\wedge y\}\hookrightarrow\{O'\})\circ\bigcup \mathrm{Im}(S_O(x),S_{O\wedge y}(y)):\overline{F}(x)\to\overline{F}(y)$, where $O$ denotes coverings of $x$ and $O'$ denotes coverings of $y$, note that here we are using Proposition 5.4.2 and 5.4.3. Easy to know that $\overline{F}$ is a presheaf and $\theta_F:F\to\overline{F}$ is a nat.

Now we prove that $\overline{F}$ is a sheaf.

The method to construct $\overline{\alpha}$ is very similar to that of $\theta_F$, and we don't bother to write it down.
\end{proof}
\paragraph{}$\overline{F}$ is called the sheafification/sheafing of $F$ or the sheaf associated to $F$, it's a process which removes some sections from and adds some sections to presheaves on the basis of the gluing axiom. Notice that, if $\alpha:F\to G$ is a presheaf morphism, then $\theta_G\circ\alpha:F\to\overline{G}$, so there is the unique $\overline{\alpha}:=\overline{\theta_G\circ\alpha}:\overline{F}\to\overline{G}$, so we turn the sheafification operation into a functor $\overline{-}:PSh^\mathcal{C}(X)\to Sh^\mathcal{C}(X)$. 
\paragraph{Proposition 5.4.5}Let $f:Y\to X$ a topology algebra homomorphism, $F:X\to\mathcal{C}$ a presheaf, then $\overline{F}\circ f=\overline{F\circ f}$. This provides a functor 
\[\overline{-}:PSh^\mathcal{C}\to Sh^\mathcal{C}:=\begin{cases}(X,F)\mapsto (X,\overline{F})\\(f:Y\to X,\alpha:F\circ f\to G)\mapsto(f,\overline{\alpha}:\overline{F\circ f}\to\overline{G})\end{cases}.\]
\paragraph{}We can discuss quotient sheaves on the basis of sheafification. If $\mathcal{C}$ is a concrete category with quotient, which is usually $Set$, $Ab$, $Rng$ and $Cat/CAT$, let $F:X\to\mathcal{C}$ a sheaf, and $F/\sim$ is one of $F$'s quotient presheaves, that is, there is a natural projection $p:F\to F/\sim$, where each $p(x):F(x)\to F(x)/\sim=(F/\sim)(x)$ is a quotient projection, then we can naturally define the quotient sheaf as the sheafification $\overline{F/\sim}$.

\section{Topology and Geometry}

\subsection{Topological Spaces}

\paragraph{Definition 6.1.1}Let $X\in TAlg$. Then topological spaces on $X$ is apex-$Set$-cosheaved algebras, and continuous mappings are $Set$-cosheaved algebra homomorphisms.
\paragraph{Theorem 6.1.1}Topology spaces defined in Definition 6.1.1 (we call it the sheaf-theoretical definition) is equivalent to the classical definition (or set-theoretical definition).
\begin{proof}
$\Leftarrow$ 

\textbf{Objects.} Let $(M,OP_M)$ is a topological space in the sense of set theory, where $OP_M$ is the family of open sets, it's a topology algebra, meanwhile an inclusion functor $in_M:OP_M\to Set$, which is actually an apex cosheaf:

If $\{f_\alpha:U_\alpha\to A\}$ making $f_\alpha\circ(U_\alpha\cap U_\beta\hookrightarrow U_\alpha)=f_\beta\circ (U_\alpha\cap U_\beta\hookrightarrow U_\beta)$, then it means $\forall a\in U_\alpha\cap U_\beta,\;f_\alpha(a)=f_\beta(a)$, namely they are compatible to each other: $\forall\alpha,\beta,\;f_\alpha|_{U_\alpha\cap U_\beta}=f_\beta|_{U_\alpha\cap U_\beta}$, so they can be glued to a larger function $f:\bigcup U_\alpha\to A:=a\mapsto f_\alpha(a)$ if $a\in U_\alpha$, easy to know that it's the only fucntion making that $\forall\alpha,\;f_\alpha=f\circ(U_\alpha\hookrightarrow \bigcup U_\alpha)$, therefore we get $(U_\alpha\hookrightarrow \bigcup U_\alpha)=paired-pushout(U_\alpha\wedge U_\beta\hookrightarrow U_\alpha,U_\alpha\cap U_\beta\hookrightarrow U_\beta)$, which is exactly the gluing axiom.

The minimum of $OP_M$ is the empty set, which is exactly the initial object of $Set$, now we know $in_M$ is a cosheaf. 

Inclusion functions $in_M(U\subseteq V)\equiv U\hookrightarrow V$ are injections, namely monos in $Set$; if there is an iso, namely bijection, between open sets $f:U\cong V$ such that $(V\hookrightarrow M)\circ f=U\hookrightarrow M$, then images of $V\hookrightarrow M$ and $U\hookrightarrow M$ are the same, namely $U=V$. So $in_M$ is apex.

%Let $U,V\subseteq W$, if $f:A\to U$ and $g:A\to V$ making $i_{U,W}\circ f=i_{V,W}\circ g$, then it means $\forall a\in A,f(a)=g(a)$, togethering with $f(a)\in U,g(a)\in V$, we get $f(a)=g(a)\in U\cap V\equiv U\wedge V$, so there exists a unique $\varphi:A\to U\wedge V$ making $f=i_{U\wedge V,U}\circ\varphi$ and $g=i_{U\wedge V,V}\circ\varphi$, therefore $(U\wedge V\hookrightarrow U,U\wedge V\hookrightarrow U)=pullback(U\hookrightarrow W,V\hookrightarrow W)$; 

\textbf{Morphisms.} Let $f:(M,OP_M)\to(N,OP_N)$ a continuous mapping in the sense of set theory, then easy to see that $f^{*}:OP_N\to OP_M:=U\mapsto f^{-1}(U)$ is an algebra homonorphism, and $\alpha_f:in_M\circ f^{*}\to in_N:=U\mapsto f|_{f^{*}(U)}$ is a nat, so $(f^{*},\alpha_f)$ is a cosheaved algebra homomorphsim, namely a continuous mapping in the sense of sheaf theory.\\

$\Rightarrow$

\textbf{Objects.} Let $(X,F)$ a topological space in the sense of sheaf theory, we write $F(x\leqslant y)$ as $i_{x,y}$ for convenience. Consider the family of subsets $OP_{F(1_X)}:=\{\mathrm{Im}\;i_{x,1_X}|x\in X\}$ of $F(1_X)$, it's actually a topology structure:

According to the gluing axiom $(i_{x_\alpha,x=\bigvee x_\alpha})=paired-pushout(i_{x_\alpha\wedge x_\beta,x_\alpha},i_{x_\alpha\wedge x_\beta,x_\beta})$, the first picture of colimits (Theorem 2.2), and the depiction of coequalizers on $Set$ (Lemma 4.1.1), we have $F(x)=\coprod F(x_\alpha)/\sim_R$, the coequalizer as a projection $e:\coprod F(x_\alpha)\to F(x)$, and embedding mappings $\{i_\alpha:F(x_\alpha)\to\coprod F(x_\alpha)\}$. Now we make a series of computation: $\bigcup \mathrm{Im}\;i_{x_\alpha,1_X}=\bigcup \mathrm{Im}\;(i_{x,1_X}\circ e\circ i_\alpha)=\bigcup i_{x,1_X}(e(\mathrm{Im}\;i_\alpha))=i_{x,1_X}(e(\bigcup\mathrm{Im}\;i_\alpha))=i_{x,1_X}(e(\coprod F(x_\alpha)))=i_{x,1_X}(F(x))=\mathrm{Im}\;i_{x,1_X}$, so we actaully obtain one topology axiom $\bigcup \mathrm{Im}\;i_{x_\alpha,1_X}=\mathrm{Im}\;i_{\bigvee x_\alpha,1_X}$.

From the last paragraph we know $F(x\vee y)=F(x)\sqcup F(y)/\sim_R$, where $aRb$ is defined as there exists $c\in F(x\wedge y)$ making that $(i_1\circ i_{x\wedge y,x})(c)=a$ and $(i_2\circ i_{x\wedge y,y})(c)=b$. Note that $i_{x\wedge y,x}$ and $i_{x\wedge y,y}$ are injective, and the embedding mappings $i_1:F(x)\to F(x)\sqcup F(x)$ and $i_2:F(y)\to F(x)\sqcup F(x)$ are injective, so $\sim_R$ is just $R$ itself, namely $F(x\vee y)=F(x)\sqcup F(y)/(i_{x\wedge y,x}(c)=i_{x\wedge y,y}(c))$. Now we make a series of computation: $\mathrm{Im}\;i_{x,1_X}\cap\mathrm{Im}\;i_{y,1_X}=\mathrm{Im}(i_{x,1_X}\circ e\circ i_1)\cap\mathrm{Im}(i_{y,1_X}\circ e\circ i_2)=i_{x,1_X}(\mathrm{Im}\;e\circ i_1)\cap i_{y,1_X}(\mathrm{Im}\;e\circ i_2)$, because $i_{x\vee y,1_X}$ is injective, so the front$=i_{x\vee y,1_X}(\mathrm{Im}\;e\circ i_1\cap\mathrm{Im}\;e\circ i_2)=i_{x\vee y,1_X}(\mathrm{Im}\;i_{x\wedge y,x\vee y})=\mathrm{Im}\;i_{x\wedge y,1_X}$, so we actaully obtain another topology axiom $\mathrm{Im}\;i_{x,1_X}\cap\mathrm{Im}\;i_{y,1_X}=\mathrm{Im}\;i_{x\wedge y,1_X}$.

Note that in the above process of proving, we use these two results: $\bigcup f(A_\alpha)=f(\bigcup A_\alpha)$, and if $f$ is injective, then $\bigcap f(A_\alpha)=f(\bigcap A_\alpha)$.

\textbf{Morphisms.} Let $(f,\alpha):(X,F)\to(Y,G)$ a cosheaved algebra homomorphism, then $\alpha(1_Y):F(1_X)\to G(1_Y)$ is actually a continuous mapping: $\alpha(1_Y)^{-1}(\mathrm{Im}\;G(y\leqslant 1_Y))=\mathrm{Im}\;F(f(y)\leqslant 1_X)$.\\

Draw a picture to make these two processes clear:
\[\begin{tikzpicture}
\path(0.5,1)node(a){Set-theoretical spaces and mappings}(0.5,-1)node(b){Sheaf-theoretical spaces and mappings}(6,1)node(c){$(M,OP_M)$}(6,-1)node(d){$(OP_M,in_M)$}(12,1)node(e){$(F(1_X),OP_{F(1_X)})$}(12,-1)node(f){$(X,F)$}
(8,1)node(g){$f$}(8,-1)node(h){$(f^*,\alpha_f)$}(14.5,1)node(i){$\alpha(1_Y)$}(14.5,-1)node(j){$(f,\alpha)$};
\draw[->](c)--(d);
\draw[->](f)--(e);
\draw[->](g)--(h);
\draw[->](j)--(i);
\end{tikzpicture}\] 

Now we prove these two processes are mutually inverse in the sense of isomorphism. 

The maximum of $OP_M$ is just $M$, and $in_M(M)$ is just $M$ itself, for any open set $U\in OP_M$, $\mathrm{Im}\;in_M(U\hookrightarrow M)$ is just $U$ itself, so we restore $(OP_M,in_M)$ back to $(M,OP_M)$. 

According to conclusion we get before, togethering with all $\mathrm{Im}\;i_{x,1_X}$ are distinct due to $F$ is apex, $OP_{F(1_X)}$ is isomorphic to $X$. Since each $i_{x,1_X}$ is an injection, $F(x)$ is isomorphic to $\mathrm{Im}\;i_{x,1_X}$, therefore $(X,F)\cong(OP_{F(1_X)},in_{F(1_X)})$. 
\end{proof}

\paragraph{}In the proof of Theoren 6.1.1, we actually construct two functors between the category of topological spaces $Top$ and the category of apex-$Set$-cosheaved algebras $ACoSh^{Set}$, and they are inverse in the sense of isomorphicness. Now we want to transfer some concepts on $Top$ to $ACoSh^{Set}$.
\paragraph{Definition 6.1.2}Let $(X,F)$ and $(Y,G)$ topological spaces, $(f,\alpha):(X,F)\to(Y,G)$. Then we call $(f,\alpha)$ a quotient mapping if $f$ is an embedding functor, and $T_x$ is an inverse image of $T_f(1_Y):Patl_X\to Patl_Y$ iff $x\in f(Y)$. Further, we call $(f,\alpha)$ an absolute quotient mapping if $f=id_X$ and $\alpha$ is a retraction.
\paragraph{Proposition} Absolute quotient is quotient.
\paragraph{Definiton 6.1.3}An algebra $X$ is called seperatable, if for any $p\ne q\in Patl_X$, there exist $x\in p$ and $y\in q$ such that $x\wedge y=0_X$. A topological space $(M,OP_M)$ or $(X,F)$ is called seperatable, if $OP_M$ or $X$ is seperatable.
\paragraph{Lemma 6.1.1}Let $X$ a seperatable algebra, then $\lim\limits _{\longleftarrow}T^X|_p=\Bigl(\{p\},(x\in p)\mapsto(\{p\}\subseteq T_x)\Bigr)$, namely any costalk $(T^X)_p=\{p\}$ is a single point set.
\begin{proof}
$T^X(x)$ are just some sets with inclusion mappings, so easy to see that $\lim\limits _{\longleftarrow}T^X|_p=\bigcap\limits _{x\in p} F(x)$, which at least has an element $p$, suppose it has another element $q$, which means for all $x\in p$, $q\in T_x$, namely $x\in q$, so $p\subset q$, therefore for any $x\in p$ and $y\in q$, it  has $x\wedge y\in q$, which makes a contradiction with seperatablity.  
\end{proof}
\paragraph{}Spaces whose all costalks are single point sets are called thin. 
\paragraph{Lemma 6.1.2}Let $(X,F)$ a seperatable topological space, and for $p\in Patl_X$, $\lim\limits _{\longrightarrow}F|_p=(F_p,\delta_p)$. Then for all $x\in X$, $\{\mathrm{Im}\;\delta_p(x)\}_{p\in T_x}$ is a partition of $F(x)$.
\begin{proof}According to the proof of Theorem 6.1.1, $(X,F)\cong(OP_{F(1_X)},in_{F(1_X)})$, so we can consider $(OP_{F(1_X)},in_{F(1_X)})$ instead of $(X,F)$ to get the same result. Because $p\in Patl(Top_F(1_X))$ is some open sets with inclusion mappings, so easy to see that $(in_F(1_X))_p=\bigcap\limits _{U\in p}U$, and $\delta_p(U)$ is the inclusion mapping which means $\mathrm{Im}\;\delta_p(U)$ is $(in_F(1_X))_p=\bigcap\limits _{U\in p}U$ itself. Firstly, for any point $a\in F(1_X)$, there is a particle $p_a=\{U|a\in U\}$, such that $a\in \bigcap\limits _{U\in p_a}U$. Secondly, for any two unequal particles $p$ and $q$, $(\bigcap\limits _{U\in p}U)\bigcap(\bigcap\limits _{U\in q}U)=\bigcap\limits _{U\in p\cup q}U=\emptyset$, because there are $U\in p$ and $V\in q$ such that $U\cap V=\emptyset$.
\end{proof}
\paragraph{Theorem 6.1.2}Let $(X,F)$ a seperatable space, then there is an absolute quotient mapping $p:(X,F)\to(X,T^X)$.
\begin{proof}Let $p\in Patl_X$. For each $x\in X$, we define a function $\alpha(x):F(x)\to T_x:=(a\in\mathrm{Im}\;\delta_p(x))\mapsto p$, according to Lemma 6.1.1, it's rational, and easy to see that $\alpha:F\to T^X$ is a nat. For each $p\in Patl_X$, we choose a representative elemtent $a_p\in \mathrm{Im}\;\delta_p(1_X)$, notice that if $p\in T_x$, $\mathrm{Im}\;\delta_p(x)\subseteq \mathrm{Im}\;\delta_p(1_X)$, so we can define a function $\beta(x):T_x\to F(x):=p\mapsto a_p$, which is a right inverse of $\alpha(x)$, easy to see that $\beta$ is a nat and is a right inverse of $\alpha$, so $(id_X,\alpha)$ is an absolute quotient mapping. 
\end{proof}
\paragraph{Lemma 6.1.3}If $(M,OP_M)$ is seperatable and sober, then it's thin.
\begin{proof}All particles have the form $p_a$, so its costalk on $p_a$ is some set containing $a$, if this set contains another point $b$, easy to see that $p_a\subset p_b$ due to $P_{-}$ is injective, but it's impossible because $M$ is seperatable, therefore the costalk on $p_a$ is $\{a\}$.
\end{proof}
\paragraph{Lemma 6.1.4}If $(M,OP_M)$ is seperatable and thin, then it's sober.
\begin{proof}For a particle $p$, if its costalk is $\{a\}$, then $p\subseteq p_a$, and due to seperatablity, $p=p_a$. If $p_a=p_b$, then $\{a\}=\{b\}$, namely $a=b$.
\end{proof}
\paragraph{Theorem 6.1.3}A topological space $(M,OP_M)$ is Hausdorff iff seperatable and sober.
\begin{proof}
$\Rightarrow$ According to Proposition 5.3.3 and 5.3.2, Hausdroff implies sober, which implies $p_{-}$ is surjective, so all $OP_M$'s particles have the form $p_a$, then $OP_M$ is seperated due to it's Hausdorff.

$\Leftarrow$ Because it's sober, so for any distinct $a\ne b\in M$, $p_a$ and $p_b$ are distinct particles of $OP_M$. because $OP_M$ is seperatable, so there are $U\in p_a$ and $V\in p_b$, namely $a\in U$ and $b\in V$, such that $U\cap V=\emptyset$. 
\end{proof}
\paragraph{Theorem 6.1.3'}A topological space $(M,OP_M)$ is Hausdorff iff seperatable and thin.
\begin{proof}According to Lemma 6.1.3 and Lemma 7.14, we know that sober and seperatable iff thin and seperatable, now we easily see it from Theorem 6.1.3.
\end{proof}
\paragraph{}These Theorems directly point out that Hausdorff spaces and seperatable topology algebras are equivalent concepts, so now we can define a sheaf-theoretical space $(X,F)$ is Hausdorff as seperatable and thin, or isomorphic to $(X,T^X)$.  

\subsection{$\mathcal{C}$-Manifolds}
\paragraph{}In the last section 8.1, we rewrite the definition of topological spaces using sheaf-theoretical language; in the chapter 5, we define a category $Smo$ which seems to help us to translate concepts about differential objects into category-theoretical language. Notice that there is an embedding functor from $Smo$ to $Set$, so $Smo$ can be seen as a subcategory of $Set$, therefore we hope to define differetial manifolds in sheaf theory by imitating topological spaces. The most simple idea is to define smooth manifolds as apex-$Smo$-cosheaved algebras, they actually are smooth manifolds but merely open sets of Euclidean spaces.
\paragraph{}Let $(M,OP_M)$ be a topological spaces. An atlas of $M$ is a family of compatible charts $\{(U_\alpha\in OP_M,\;F_\alpha:U_\alpha\tilde{\to} V_\alpha\in OP_{\mathbb{R}^n})\}_{\alpha\in J}$, where each $F_\alpha$ is a homeomorphism and $F_\alpha\circ F_\beta|_{U_\alpha\cap U_\beta}$ is smooth, if we regard each $F_\alpha$ as a local apex-$Smo$-cosheaf on $U_\alpha^\leqslant$ whose $F_\alpha(U_\alpha)=V_\alpha$, we may be able to use the method of gluing local cosheaves to construct smooth structure. So we make a try:
\paragraph{}Let $M=(X,F)$ a topological space. $\{U_\alpha\}_{\alpha\in J}$ is a covering of $1_X$ and each $U_\alpha^\leqslant$ equips with an apex-$Smo_n$-cosheaf $F_\alpha$. For each $\alpha\in J$, there is a natural isomorphism $\varphi_\alpha:F|_{U_\alpha}\to F_\alpha$, so for each $\alpha,\beta\in J$, we can build a natural isomorphism (\textbf{equivalentor}) $\varphi_{\alpha\beta}=(\varphi_\beta\circ\varphi^{-1}_\alpha)|_{U_\alpha\wedge U_\beta}:F_\alpha |_{U_\alpha\wedge U_\beta}\to F_\beta |_{U_\alpha\wedge U_\beta}$, and they satisfy \textbf{the cochain condition with restrictions}
\[\begin{cases}\varphi_{\alpha\alpha}=id_{F_\alpha},\\
\varphi_{\beta\alpha}\circ\varphi_{\alpha\beta}=id_{F_\alpha|_{U_\alpha\wedge U_\beta}},\\
\varphi_{\gamma\alpha}\circ\varphi_{\beta\gamma}\circ\varphi_{\alpha\beta}|_{U_\alpha\wedge U_\beta\wedge U_\gamma}=id_{F_\alpha|_{U_\alpha\wedge U_\beta\wedge U_\gamma}},\end{cases}\]
 The other way round, given  a family of charts $\{(U_\alpha,\;F_\alpha:U_\alpha\to Smo)\}_{\alpha\in J}$ and a cochain group $\{\varphi_{\alpha\beta}:F_\alpha|_{U_\alpha\wedge U_\beta}\to F_\beta|_{U_\alpha\wedge U_\beta}\}_{\alpha,\beta\in J}$, we can define a topological space: firstly we adopt a choice function $ch:\mathrm{Ob}X\cup\mathrm{Mor}X\to J$ making that $V\in \mathrm{Ob}\;U^\leqslant_{ch(V)}$ and $f\in\mathrm{Mor}\;U^\leqslant_{ch(f)}$, then we define $PF(V):=F_{ch(V)}(V)$ and $PF(f=V\leqslant  U):=\varphi_{ch(f),ch(U)}(U)\circ F_{ch(f)}(f)\circ\varphi_{ch(V),ch(f)}(V)$, easy to prove that $PF:\bigcup_{\alpha\in J} U^\leqslant_\alpha\to Smo$ is a functor, then according to Theorem 5.2.1 (gluing of local cosheaves), we get a cosheaf $F:X\to Set$, and natural isomorphisms $\varphi_\alpha:F|_{U_\alpha}\to F_\alpha:=V\mapsto\varphi_{ch(V),\alpha}(V)$ which satisfy $\varphi_{\alpha\beta}=(\varphi_\beta\circ\varphi^{-1}_\alpha)|_{U_\alpha\wedge U_\beta}$, easy to prove that such $F$ is unique in the sense of isomorphicness. Now we see that the cochain condition is the connotation of equivalentors. Therefore, we can regard a smooth manifold as a topological space which has a family of apex-$Smo$-cosheaves equipped with a cochain group as a covering (in the sense of isomorphicness). This definition is merely an indiscriminately copy of the classical definition of smooth manifolds, so we may be unable to discover new information from it, we seek for an expression instead which regards all manifolds as a holistic mathematical object and is more correlated to category theory.
\paragraph{}There is another way to construct a cosheaf which may represent smooth manifolds. Firstly we add a requirement: images of all $F_\alpha$ are disjoint. For any $V\in\bigcup_{\alpha\in J}U_\alpha^\leqslant$, there is a cochain group $\{\varphi_{\alpha\beta}(V)|V\leqslant U_\alpha\wedge U_\beta\}$, which has declared an equivalence class $\{F_\alpha(V)|V\leqslant U_\alpha\}$. Therefore we can get a sketch quotient category of $Smo$ and the quotient functor $\pi:Smo\to Smo/\sim$, easy to know that $\{\pi\circ F_\alpha:U_\alpha\to Smo/\sim\hookrightarrow Set/\sim\}_{\alpha\in J}$ is a family of compatible cosheaves, because $Set$ is complete, according to Theorem 5.2.1, they can be glued into a cosheaf $F:X\to Set/\sim$. Construction of morphisms is a little complicated: let $(X,F:X\to Set/\sim_F)$ and $(Y,G:X\to Set/\sim_G)$ be smooth manifolds where images of $\{F_\alpha\}$ and $\{G_\alpha\}$ are disjoint, notice that there naturally are two quotient functors $p_F:Set/\sim_F\to Set/\sim_F/\sim_G$ and  $p_G:Set/\sim_F\to Set/\sim_F/\sim_G$. Then a smooth mapping from $(X,F)$ to $(Y,G)$ is a sheaved algebra homomorphism from $(X,p_F\circ F)$ to $(Y,p_G\circ G)$. This definition has many drawbacks: 1. It relies on the harsh condition that images of $F_\alpha$ are disjoint; 2. It relies on solving quotient categories many times, we even need to solve them three and four times when we talk about the composition of smooth mappings and its associative law; 3. It relies on a concrete atlas and doesn't indicate which atlases are equivalent, namely representing the same smooth structure. In a word, it depends on the method of concretely gluing.
\paragraph{}But the method above of gluing local cosheaves enlightens us to introduce a concept of abstractly gluing by regarding local cosheaves as local sections of an another sheaf, and the abstractly gluing is provided by the gluing axiom. To discuss further, we firstly rewrite Theorem 5.2.1 (gluing of local sheaves):
\paragraph{Theoren 6.2.1}Let $\mathcal{C}$ a paired-pushout-complete category and has the initial object.
\subparagraph{$\bf{1.}^\mathrm{Ob}$}Let $X\in TAlg$, then $CAT$-presheaf
\[Sh_{lg}^\mathcal{C}(X):X\to CAT:=\begin{cases}x\mapsto Sh^\mathcal{C}(x^\leqslant)\\x\leqslant y\mapsto\Biggl(\begin{cases}F\mapsto F|_{x^\leqslant}\\\alpha\mapsto\alpha|_{x^\leqslant}\end{cases}:Sh^\mathcal{C}(y^\leqslant)\to Sh^\mathcal{C}(x^\leqslant)\Biggr)\end{cases} \]
is a sheaf. Easy to see that a $\mathcal{C}$-sheaf on $X$ is a  global object-section $F\in\mathrm{Ob}(Sh_{lg}^\mathcal{C}(X)(1_X))$.
\subparagraph{$\bf{2.}^\mathrm{Mor}$}Let $f\in Hom_{TAlg}(Y,X)$, then $(f,Sh_{lg}^\mathcal{C}(f):Sh_{lg}^\mathcal{C}(X)\circ f\to Sh_{lg}^\mathcal{C}(Y))$ is a sheaved algebra homomorphism, where 
\[Sh_{lg}^\mathcal{C}(f):=y\mapsto\begin{cases}F\mapsto F\circ(f|_{y^\leqslant})\\\alpha\mapsto\alpha\circ(f|_{y^\leqslant})\end{cases} \]
is a nat. Easy to see that if $(f,\alpha:F\circ f\to G)$ is a sheaved algebra homomorphism, then $F\circ f=Sh_{lg}^\mathcal{C}(f)(1_Y)(F)$ and $\alpha\in\mathrm{Mor}(Sh_{lg}^\mathcal{C}(Y)(1_Y))$ is a global morphism-section.
\\\\Replacing $Sh$ to $CoSh$ the theorem is still true. $(Co)Sh_{lg}^\mathcal{C}(X)$ is called \textbf{the category sheaf of local $\mathcal{C}$-(co)sheaves on $X$}, where "lg" means "locally glue". The theorem actually provides the contravariant functors $(Co)Sh_{lg}^\mathcal{C}(-):TAlg\to Sh^{CAT}$.
\paragraph{}If $F_1\in CoSh_{lg}^\mathcal{C}(x^\leqslant)$ and $F_2\in CoSh_{lg}^\mathcal{C}(y^\leqslant)$ are two object-sections of $CoSh_{lg}^\mathcal{C}(X)$, then $F_1$ and $F_2$ are compatible if $F_1|_{(x\wedge y)^\leqslant}=F_2|_{(x\wedge y)^\leqslant}$. But we find a subtle difference between the compatibility in Theorem 6.2.1 and the compatibility we need, which is $F_1|_{(x\wedge y)^\leqslant}\cong F_2|_{(x\wedge y)^\leqslant}$, and what is more annoying is that a family of local cosheaves compatible to each other can't always be glued to a larger cosheaf. However, we just need to change the definition of $Cosh_{lg}^\mathcal{C}(X)$ slightly to get what we need with the help of sheafication which provides the concept of abstractly gluing.
\paragraph{}Consider a presheaf in form
\[CoSh_{lg}^{\mathcal{C},p}(X):X\to CAT:=\begin{cases}x\mapsto CoSh^\mathcal{C}(x^\leqslant)/\cong\\x\leqslant y\mapsto\begin{cases}[F]\mapsto [F|_{x^\leqslant}]\\ [\alpha ]\mapsto[\alpha|_{x^\leqslant}]\end{cases}\end{cases},\]
where $\cong$ represents natural isomorphism, which means each $CoSh_{lg}^{\mathcal{C},p}(X)(x)$ is a skeleton quotient category. Obviously, isomorphic cosheaves' restrictions are isomorphic, but equivalent nats' restrictions are not necessarily equivalent, so we require that all $CoSh^\mathcal{C}(x^\leqslant)$'s choices of cochain groups must lead to $\alpha\sim^\mathrm{Mor}\beta\Rightarrow\alpha|_{x^\leqslant}\sim^\mathrm{Mor}\beta|_{x^\leqslant}$, where $\alpha,\beta:y^\leqslant\to\mathrm{Mor}(\mathcal{C})$ ara nats between $\mathcal{C}$-cosheaves on $y^\leqslant$, from now on, functors $CoSh_{lg}^{\mathcal{C},p}(X)(x\leqslant y)=\begin{cases}[F]\mapsto [F|_{x^\leqslant}]\\ [\alpha ]\mapsto[\alpha|_{x^\leqslant}]\end{cases}$ get their rationality, and we can naturally write $[F|_{x^\leqslant}]$ and $[\alpha|_{x^\leqslant}]$ as $[F]|_{x^\leqslant}$ and $[\alpha]|_{x^\leqslant}$, meanwhile there is a natural projection $p:CoSh^\mathcal{C}_{lg}(X)\to CoSh_{lg}^{\mathcal{C},p}(X)$, now we see that ways to quotient is determined by natural projections, so we write the quotient presheaf determined by $p$ as $CoSh_{lg}^{\mathcal{C},p}(X)$ for convenience.  
\paragraph{}Let $\{[F_\alpha]\in CoSh_{lg}^{\mathcal{C},p}(X)(U_\alpha)\}_{\alpha\in J}$ is a family of compatible object-sections, which means $[F_\alpha]|_{U_\alpha\wedge U_\beta}=[F_\beta]|_{U_\alpha\wedge U_\beta}$, then choose any representative elements $\{F_\alpha\}$, there naturally exists a group of natural isomorphisms $\{\varphi_{\alpha\beta}:F_\alpha|_{U_\alpha\wedge U_\beta}\to F_\beta|_{U_\alpha\wedge U_\beta}\}$ which are members in the cochain groups of $[F_\alpha|_{U_\alpha\wedge U_\beta}]$s, in other word, in the equivalence classes $[id_{F_\alpha|_{U_\alpha\wedge U_\beta}}]$s, notice that functors $CoSh_{lg}^{\mathcal{C},p}(X)(x\leqslant y)$ are identity-preservative, in other word, functors $CoSh^{\mathcal{C},p}_{lg}(X)(x\leqslant y)$ map members in cochain groups to members in cochain groups, so $\{\varphi_{\alpha\beta}\}$ satisfies the cochain condition with restrictions, which means $\{F_\alpha\}$ and $\{\varphi_{\alpha\beta}\}$ gives a smooth structure. Now we explain that all choices of family of compatible object-sections represent the same smooth structure: let $\{F_\alpha\}_{\alpha\in J}$ and $\{G_\alpha\}_{\alpha\in J}$ are two choices, $\{\varphi_{\alpha\beta}\}$ and $\{\psi_{\alpha\beta}\}$ are their cochain groups, and $\{\omega_\alpha:F_\alpha\to G_\alpha\in[id_{F_\alpha}]\}$ are isomorphisms, easy to see that $\psi_{\alpha\beta}=\omega_\beta|_{U_\alpha\wedge U_\beta}\circ\varphi_{\alpha\beta}\circ(\omega_\alpha|_{U_\alpha\wedge U_\beta})^{-1}$, in reverse, $\varphi_{\alpha\beta}=(\omega_\beta|_{U_\alpha\wedge U_\beta})^{-1}\circ\psi_{\alpha\beta}\circ\omega_\alpha|_{U_\alpha\wedge U_\beta}$, so they can be transformed to each other by some isomorphisms, which means they are  essentially the same. Therefore, we obtain:
\paragraph{Definition $\bf{6.2.1}^\mathrm{Ob}$}Let $X\in TAlg$, then a quotient presheaf $CoSh_{lg}^{\mathcal{C},p}(X)$'s sheafication $\overline{CoSh_{lg}^{\mathcal{C},p}(X)}$ is called \textbf{a category sheaf of local $\mathcal{C}$-manifolds on $X$}. Its global object-section $\Sigma\in\mathrm{Ob}\;\overline{CoSh_{lg}^{\mathcal{C},p}(X)}(1_X)$ is called  a $\mathcal{C}$-manifold structure on $X$, and two-tuples $(X,\Sigma)$ a $\mathcal{C}$-manifold. And we call a covering $\{(x_\alpha,\Sigma|_{x_\alpha})\}$ of $(X,\Sigma)$ an atlas of it, if all $\Sigma_\alpha$ are local cosheaves.
\paragraph{}Now we see that a category sheaf of local manifolds is namely a skeleton quotient sheaf of $CoSh_{lg}$. The $\mathcal{C}$-manifolds defined here are the most general situation, of course, we can put some restrictions on them, such as apex, Hausdorff and so on to reach the demand of various concrete mathematical researches.
\paragraph{}To discuss morphisms between $\mathcal{C}$-manifolds, we need promote $CoSh_{lg}^\mathcal{C}(f)$ to $\overline{CoSh_{lg}^{\mathcal{C},p,q}(f)}$. Firstly easy to construct the functor
\[CoSh_{lg}^{\mathcal{C},p,q}(f):CoSh_{lg}^{\mathcal{C},p}(X)\circ f\to CoSh_{lg}^{\mathcal{C},q}(Y):=y\mapsto\begin{cases}[F]\mapsto[F\circ(f|_{y^\leqslant})]\\ [\alpha]\mapsto[\alpha\circ(f|_{y^\leqslant})]\end{cases},\]
if $p$ and $q$ satisfy $\forall y\in Y,\;\alpha\sim^\mathrm{Mor}\beta\Rightarrow\alpha\circ(f|_{y\leqslant})\sim^\mathrm{Mor}\beta\circ(f|_{y\leqslant})$. Then using the functor $\overline{-}:PSh^{CAT}\to Sh^{CAT}$ we have $\overline{CoSh_{lg}^{\mathcal{C},p,q}(f)}:\overline{CoSh_{lg}^{\mathcal{C},p}(X)\circ f}\to\overline{CoSh_{lg}^{\mathcal{C},q}(Y)}$, where $\overline{CoSh_{lg}^{\mathcal{C},p}(X)\circ f}$ is equal to $\overline{CoSh_{lg}^{\mathcal{C},p}(X)}\circ f$.
\paragraph{Definition $\bf{6.2.1}^\mathrm{Mor}$}Let $(X,\Sigma)$ and $(Y,\Gamma)$ $\mathcal{C}$-manifolds. Then a morphism from $(X,\Sigma)$ to $(Y,\Gamma)$ is a two-tuples $\bigl(f:Y\to X,\;\alpha:\overline{CoSh_{lg}^{\mathcal{C},p,q}(f)}(1_Y)(\Sigma)\to\Gamma\bigr)$, namely $\alpha$ is a global morphism-section $\alpha\in \mathrm{Mor}\;\overline{CoSh_{lg}^{\mathcal{C},q}(Y)}(1_Y)$.
\\\\Then we get the category of $\mathcal{C}$-manifolds $Man^\mathcal{C}$. Particularly, if $\mathcal{C}$ is a subcategory of $Set$, then $\mathcal{C}$-manifolds can be regarded as topological spaces: $CoSh_{lg}^{Set}(X)/\cong$ is already a sheaf, that is, local $Set$-cosheaves always can be glued to a larger cosheaf. This characteristic depends on some sorts of inner symmetry of $Set$.
\paragraph{}Let $(X,\Sigma)\in Man^\mathcal{C}$, there is an embedding functor \[in_\Sigma:X\to Man^\mathcal{C}:=\begin{cases}x\mapsto (Sub(x),\Sigma|_x)\\x\leqslant y\mapsto(Sub(x\leqslant y),z\in y^\leqslant\mapsto \Sigma(z\wedge x\leqslant z)) \end{cases},\]which is a cosheaf.

\subsection{Smooth Manifolds}

\paragraph{}A smooth manifold is an apex-$Smo$-manifold, that is, any its atlas is a family of apex cosheaves. Of course, we can put on some other conditions to reach concrete researches, such as compact, Hausdorff and having a countable basis. Subcategories $Smo_n$ of $Smo$ don't have isomorphic objects mutually, so an atlas is always a family of $Smo_n$-cosheaves with a specific index n, and we call it $Smo$-manifold's dimension.
\paragraph{Lemma 6.3.1}Let $F:X\to\mathcal{C}$ a cosheaf, then for any $A\in\mathcal{C}$, $Hom^A(-)\circ F$ is a $Set$-sheaf. Dually, let $F:X\to\mathcal{C}$ a sheaf, then for any $A\in\mathcal{C}$, $Hom_A(-)\circ F$ is a $Set$-sheaf.
\begin{proof}We only prove the first and the second is similar. The gluing axiom $Hom^A(-)\circ F$ on is a direct corollary of Theorem 4.2.3. Because $F(0_X)$ is the initial object, so $Hom(F(0_X),A)$ is a single point set, namely the terminal object of $Set$.
\end{proof}
\paragraph{}If $A$ is a $\mathcal{C}$-object, where $\mathcal{C}$ is a concrete category whose terminal object appears as a single point set, then the lemma is still true, especially $Ab$, $Rng$ and $Cat/CAT$.
\paragraph{}There are special objects on $Man^{Smo}$. The 1-dimensional Euclidean space $\mathbb{R}^1$ is a ring object on $Man^{Smo}$, so if $M=(X,\sigma)\in Man^{Smo}$, then $Hom^\mathbb{R}(-)\circ in_\sigma$ is the sheaf of smooth scalar fields $F(M)$.

\end{document}